\documentclass[preprint,3p]{elsarticle}

\biboptions{sort&compress}

\usepackage{amsfonts}
\usepackage{amsmath}
\usepackage{amssymb}
\usepackage{amsthm}
\usepackage{graphicx}
\usepackage{tabularx}
\usepackage{cleveref}
\usepackage{color}

\newcommand{\x}{\mathbf{x}}

\newcommand{\N}{\mathcal{N}}

\newcommand{\p}{\partial}
\newcommand{\yr}{y^{\mathrm{ref}}}

\newcommand{\EE}{\mathbb{E}}
\newcommand{\PP}{\mathbb{P}}
\newcommand{\RR}{\mathbb{R}}

\newcommand{\<}{\langle}
\renewcommand{\>}{\rangle}

\journal{Physica D}

\begin{document}

\begin{frontmatter}



\title{Controlled fluid transport by the collective motion of microrotors}

\author[wisc,clem]{Jake Buzhardt}
\ead{buzhardt@wisc.edu}
\author[clem]{Phanindra Tallapragada}
\ead{ptallap@clemson.edu}

\affiliation[wisc]{organization={Department of Chemical \& Biological Engineering, University of Wisconsin-Madison},
            city={Madison}, 
            state={WI},
            country={USA}}
\affiliation[clem]{organization={Department of Mechanical Engineering, Clemson University},
            city={Clemson}, 
            state={SC},
            country={USA}}

\begin{abstract}
Torque-driven microscale swimming robots, or microrotors, hold significant potential in biomedical applications such as targeted drug delivery, minimally invasive surgery, and micromanipulation. 
This paper addresses the challenge of controlling the transport of fluid volumes using the flow fields generated by interacting groups of microrotors. 
Our approach uses polynomial chaos expansions to model the time evolution of fluid particle distributions and formulate an optimal control problem, which we solve numerically.  
We implement this framework in simulation to achieve the controlled transport of an initial fluid particle distribution to a target destination while minimizing undesirable effects such as stretching and mixing.
We consider the case where translational velocities of the rotors are directly controlled, as well as the case where only torques are controlled and the rotors move in response to the collective flow fields they generate. 
We analyze the solution of this optimal control problem by computing the Lagrangian coherent structures of the associated flow field, which 
reveal the formation of transport barriers that efficiently guide particles toward their target.  This analysis provides insights into the underlying mechanisms of controlled transport.

\end{abstract}




\end{frontmatter}

\section{Introduction}

Recent advances in microscale swimming robots have sparked a renewed interest in understanding, modeling, and controlling fluid motion in microfluidic environments \cite{Nelson2010,Sitti2015}.
Such swimming robots have been proposed as a means of transporting cargo through complex environments, which could prove useful for many biomedical applications, such as targeted therapeutics and drug delivery
, minimally invasive surgery, and other micro-manipulation tasks \cite{Erkoc2019,Xu2020, Liu2022,Chen2023}.  
While many proposed methods in this area involve physically attaching the cargo to the carrier \cite{Tottori2012,Steager2013,Ren2019}, an alternative approach is to achieve transport by using the flow field produced by a microswimmer or microrobot \cite{Pawashe2012,Petit2012,Ye2012,Zhang2012,Tung2013,Huang2014,Ye2014,Zhou2017,Lin2018,Khalil2020,Buzhardt2021}.  
In this way, microswimmers can steer volumes of fluid or submerged cargo to a target destination. 
However, due to the small size and surface areas of microrobots, many complex manipulation tasks require groups or even swarms of microrobots \cite{Yang2021}.  
In order to carry out these tasks effectively, it is necessary to develop efficient multiagent control strategies which account for the complex interactions of these swimmers with their environment including other nearby swimmers, boundaries, other cargo and passive particles, and the fluid medium itself.
This work contributes to addressing this challenge by developing a method to achieve controlled transport of fluid particle distributions using the fluid velocity fields generated by groups of mobile microrotors by controlling the torques and translational velocities of the rotors. 

In order to approach this problem from a dynamical systems perspective, we make use of a simplified model of a torque-driven microswimmer, or microrotor, where the microrotor is represented as a point-torque in a Stokes flow. 
Such microrotors are common in practice, as many artificial microswimmers achieve propulsion by a torque induced through externally applied magnetic fields \cite{Yang2021, Han2018}.  
In many of these cases, the applied torque leads to a rotation of the swimmer body, which is designed so that such rotational motion is coupled to a net translation, allowing the robot to achieve locomotion. 
It has been shown that a point torque model serves as a good approximation for such torque-driven microrotors, as the anti-symmetric, rotational component of the velocity field dominates over symmetric dipole and higher order singularity terms \cite{Buzhardt2019,Pal2022}. Following from this, the point-torque model captures both the qualitative features of the resulting fluid velocity field and hydrodynamic couplings between nearby microrotors.
While these swimmers are effective and efficient, they often require complicated geometries to achieve the necessary rotation-translation coupling, which can make them challenging to manufacture.  

An alternative approach to achieving torque-driven locomotion from a torque induced by a magnetic field is to use simple geometries such as magnetic cylinders, spheres, or bodies composed of these (see Lyu \textit{et al.} \cite{Lyu2023} for a recent review).  For such simple body geometries, the spinning motion due to the magnetic field will not, in general, lead to locomotion for an individual swimmer in free space, but these micro-rotors can achieve collective motion by group interactions\cite{Leoni2011, Fily2012} or by interacting with boundaries \cite{Tierno2008, Sing2010}.  
Grzybowski \textit{et al.} \cite{Grzybowski2000,Grzybowski2002} studied the pattern formation of such spinning cylinders, demonstrating that they tend to spontaneously assemble into ordered structures due to a combination of magnetic and hydrodynamic interactions, and similar mechanisms were later proposed as a strategy for simple and efficient microfluidic mixing \cite{Campbell2004,Lee2009, Ballard2016}. 
Driscoll \textit{et al.} \cite{Driscoll2017} showed through both experiments and numerics that the spontaneous formation of such structures can be achieved from hydrodynamic interactions alone (without other forces arising from magnetic attraction, collisions, or sense-and-response interactions), and that these collections of micro-rotors can create strong advective flows of potential use for applications requiring controlled fluid transport or transport of passive cargo.  
These experimental results have inspired numerous numerical investigations of the emergent collective phenomena and instabilities of large groups of microrotors \cite{Yeo2015, Delmotte2017, Delmotte2017a, Sprinkle2017, BalboaUsabiaga2017, Sprinkle2020}. 
Other works have studied interactions, collective behaviour, and transport arising due to small groups of microrotors \cite{Climent2007,Leoni2011, Fily2012, Lushi2015, Tallapragada2019, Delmotte2019, Guo2024}.

It has been shown repeatedly that these microrotor systems can yield surprisingly complex dynamics, despite their apparent simplicity.  
One of the earliest works to demonstrate this was by Meleshko and Aref \cite{Meleshko1996}, who proposed the \emph{blinking rotlet} model, which consists of two fixed-in-place, periodically blinking microrotors enclosed in a circular domain, as one of the simplest systems exhibiting chaotic advection \cite{Aref2017}.
Lushi and Vlahovska \cite{Lushi2015} studied the motion of groups of interacting rotors, showing that these groups can exhibit stable structures which rotate or translate collectively depending on the rotor spin configurations.  Further, they demonstrated that groups of three or fewer rotors tend to exhibit regular motions, while groups of four or more rotors tend to produce chaotic motions, a qualitatively similar behaviour to that known of point-vortex motions. 
Delmotte \cite{Delmotte2019} studied the stability and bifurcations of periodic ``leapfrog'' trajectories of rotors interacting above a plane-wall as the relative strength of the driving torques to gravitational forces is varied. 
Guo \textit{et al.} \cite{Guo2024} analyzed the different possible trajectories of a pair of interacting rotors in a circular confinement, mapping out the parameter space and considering the effects of finite rotor size.  
Tallapragada and Sudarsanam \cite{Tallapragada2019} presented a numerical study of the mixing properties of the flows produced by a pair of interacting, rotors of constant strengths confined within a circular domain. 

While these works analyze the complex and chaotic motions arising from uncontrolled systems of rotors, there has been relatively little attention paid to these systems from a control systems perspective. Here, our focus is on developing control strategies for systems of interacting, mobile microrotors to achieve a desired fluid transport. This is of particular interest in many of the applications mentioned above, which require targeted cargo transport and delivery. Specifically, we consider the problem of steering a distribution of fluid particles from a given initial configuration to a target destination while minimizing the spread of the particles. We formulate this controlled transport problem in terms of probability densities of fluid particles, where the objective is to drive the mean of the particle distribution to a target distribution while minimizing the variance of the distribution. In previous work \cite{Buzhardt2024}, we studied a similar problem of controlled transport by rotors fixed in space, using a transfer operator-based approach to formulate and solve the controlled density transport problem. There, it was shown that a pair of fixed rotors can effectively transport an ensemble of fluid particles to a target destination, though not without significant stretching of the distribution, which is intensified by the presence of boundaries. The inevitability of such stretching in that case was likely due to the limitations of flow fields that can be produced by a pair of stationary rotors. 
A related area that has received some attention is the control of vortex motion in inviscid flows \cite{Noack2000, Vainchtein2006, Protas2008, Krishna2023a}. Most of these works, however, focus on optimizing mixing within systems of vortices rather than achieving directed transport, which is the focus of our study. 

Here we show that larger groups of mobile, interacting rotors are able to generate much richer dynamics and more complex flow fields (as compared to the fixed-rotor case), thereby enabling more effective control of the fluid particle distribution. 
However, this improved ability to steer distributions comes at the cost of making the control task somewhat more challenging, due to higher dimension and more complex system dynamics.  Specifically, the dimension of the dynamical system is higher, as the positions of the rotors must be tracked as well as the position of the relevant fluid particles.  In studying this control problem with mobile rotors, we consider two models for the motion of the rotors: one where both the rotor translational velocities and the rotor strengths are control inputs and one where only the rotor strengths are controlled, with each rotor being advected by the flow produced by the other rotors. For the case of planar rotors considered here, this yields a state dimension of $2(n_r+1)$ and a control dimension of $3n_r$ in the former case and $n_r$ in the latter case, where $n_r$ is the number of rotors present. 
This increase in dimension makes it challenging to apply the transfer operator method that was used in Ref. \cite{Buzhardt2024}, as the dimension of the finite dimensional approximation of the transfer operator needed for a good approximation typically scales with the dimension of the system \cite{Klus2016}. For this reason, such operator-based methods are not usually applied for systems of dimension higher than two or three.  
Therefore, we instead adopt the method of polynomial chaos expansions \cite{Xiu2002,Xiu2010} to model the transport of probability density by the system.  This approach gives the advantage that the higher dimensional representation of the density evolution scales with the number of stochastic quantities rather than the dimension of the system. 
For the systems of microrotors considered here, this means that the number of rotors can be increased without increasing the number of basis functions needed for a good representation of the probability density.
Background and more detail on polynomial chaos methods will be given in Sec. \ref{sec:gpc}. 

The remainder of this paper is organized as follows.  
In Sec. \ref{sec:rotlet_models}, we introduce the singularity models of the fluid flow and microrotor motion.  
In Sec. \ref{sec:gpc}, we introduce the preliminaries of the polynomial chaos expansion method, as applied to a stochastic control system, including the derivation of the deterministic surrogate system by Galerkin projection, the propagation of moments, and the formulation of the optimal control problem.   
In Sec. \ref{sec:movingrotors}, we apply this method to the problem of steering a distribution of fluid particles using mobile microrotors, where the rotor strengths and translational velocities are controlled directly.  In Sec. \ref{sec:lcs}, we analyze the flow structures associated with the optimal control solution, which reveal regions of the fluid which are entrained by the moving rotors and pulled toward the target.  In Sec. \ref{sec:torquecontrol}, we study the case where the rotor velocities are not controlled directly, but instead the rotors are advected by the flow field.  
Finally, we conclude with a summary and discussion of future directions. 

\section{Dynamics of groups of microrotors} \label{sec:rotlet_models}
The motion of the microrotors under consideration here is described by very small length and velocity scales.  In such a setting, the fluid motion is well-described by the Stokes equations 
\begin{subequations}
\begin{align}
    -\nabla p + \mu \nabla^2 \mathbf{u} &= \mathbf{0}\\
    \nabla \cdot \mathbf{u} &= 0 
\end{align}
\end{subequations}
where $\mathbf{u}$ is the fluid velocity field, $p$ is the pressure field, and $\mu$ is the fluid viscosity.  
The disturbance velocity field produced by a spinning body in this setting can be approximated by a rotlet, the singularity solution of the Stokes equation associated with a point torque, given in two dimensions by \cite{Pozrikidis1992}
\begin{equation}
    \mathbf{u}(\mathbf{x}) = -\gamma\hat{k} \times \frac{~\x - \x_R~}{r^2}
\end{equation}
where $\x_R$ is the location of the rotor, $r=\|\x-\x_R\|$, $\gamma$ is the strength of the rotor, and $\hat{k}$ is the unit vector normal to the plane containing the rotors. 
The linearity of the Stokes equations implies that the velocity field due to a collection of $n_r$ rotors can be expressed as the the superposition of each of their individual velocity fields as 
\begin{equation} \label{eq:rotlets_vel_field}
    \mathbf{u}(\x) = -\sum_{i=1}^{n_r}\left(\gamma_i\hat{k} \times \frac{~\x - \x_i~}{r_i^2}\right)
\end{equation}
where $\gamma_i$ and $\x_i$ are the strength and position of the $i$-th rotor, respectively, and $r_i = \|\x-\x_i\|$. 

In what follows we will consider two models for the motion of the rotors themselves.  In the first model, we assume that the rotor's translational velocity can be prescribed directly.  
We formulate this as a control system describing the two dimensional position of a fluid particle and each rotor.  This gives the following the following dynamical system:
\begin{equation}
    \frac{d}{dt}
    \begin{bmatrix}
    x_p\\ y_p\\ x_{r1} \\ \vdots \\ x_{rn_r} \\ y_{r1} \\ \vdots \\ y_{rn_r}
    \end{bmatrix}
    =
     \begin{bmatrix}
     \vspace{10pt}
    -\sum_{i=1}^{n_r} \left(\gamma_i\hat{k} \times \frac{~\x - \x_i~}{r_i^2}\right)
    \\ v_{x1} \\ \vdots \\ v_{xn_r} \\ v_{y1} \\ \vdots \\ v_{yn_r}
    \end{bmatrix}
\end{equation}
where the fluid particle position is $\x = (x_p,y_p)$, and the rotlet positions are $\x_i = (x_{ri},y_{ri})$ for $i = 1, \dots, n_r$, and the control input $u\in\RR^{3n_r}$ is taken to be the vector of the rotor strengths and translational velocity components; that is, $u = [\gamma_1, \dots, \gamma_{nr},v_{x1},\dots, v_{xr}, v_{y1}, \dots, v_{yr}]^T$. So, we have a system of $2(n_r+1)$ states and $3n_r$ controls. 

The assumption that the translational velocities can be controlled directly is a somewhat strong assumption, as in many of the physical systems that we are motivated by, it is not possible to decouple the rotor translational motion from the rotor strengths so cleanly. 
For instance, in the magnetically driven microswimmers
\cite{Nelson2010,Sitti2015,Yang2021, Han2018, Buzhardt2019},
the swimmer's translational motion is coupled to it's rotational motion and torques as described by a hydrodynamic mobility matrix.  In cases involving multiple swimmers, their collective motion is given by a superposition of the swimmers self-mobilities and cross mobilities which describe the hydrodynamic coupling between swimmers, as encoded by a grand mobility matrix \cite{Kim2005}. 
However, from a control perspective, this simpler model is much easier to optimize, and can thus shed light on the optimal solutions for more complex scenarios where the motion of the rotors is fully coupled. 

In the second model, we move away from the assumption that the rotor velocities can be directly controlled and instead consider the case where only the rotor strengths are controlled.  In this case, the rotor velocities are determined by the flow field produced by all other rotors (with zero self-contribution). 
For this model, the fluid particle motion is described by the same velocity field as in Eq. \ref{eq:rotlets_vel_field}, specifically, 
\begin{equation}\label{eq:lonly_rotorparticle}
    \frac{d}{dt}
    \begin{bmatrix}
        x_p\\ y_p
    \end{bmatrix}
    = -\sum_{i = 1}^{n_r} \left(\gamma_i\hat{k} \times \frac{\mathbf{x} - \mathbf{x}_{ri}}{r_i^2}\right)
\end{equation} 
and the rotor motion is given by 
\begin{equation} \label{eq:lonly_rotorrotor}
    \frac{d}{dt}
    \begin{bmatrix}
        x_{rj}\\
        y_{rj}
    \end{bmatrix}
     = 
     -\sum_{\substack{i = 1\\ i \neq j}}^{n_r} \left(\gamma_i\hat{k} \times \frac{\mathbf{x}_{rj} - \mathbf{x}_{ri}}{r_{ij}^2}\right)
\end{equation}
for each rotor enumerated by $j = 1, \dots, n_r$ where $r_{ij} = \|\mathbf{x}_j - \mathbf{x}_i\|$ is the distance between the the $i^{\text{th}}$ and $j^{\text{th}}$ rotors. 
Eqs. \ref{eq:lonly_rotorparticle} and \ref{eq:lonly_rotorrotor} define a control system of $2(n_r+1)$ states with $n_r$ control inputs. The control is taken to be the vector of rotor strengths $u = [\gamma_1, \dots, \gamma_{n_r}]^T$.

With these two models we study the problem of transporting a distribution of fluid particles from a given initial distribution to a final distribution in a fixed time.  This fluid transport is modeled as transport of a probability density which describes the distribution of fluid particles of interest. 
We formulate the time evolution and optimal control problem for this problem using polynomial chaos expansions, as described in section \ref{sec:gpc}.

\section{Polynomial chaos expansions}  \label{sec:gpc}
Polynomial chaos expansions are a popular method of quantifying and propagating uncertainties which do not depend on time, that is, stochasticity associated with constant parameters or initial conditions of a system.  
The original idea of polynomial chaos expansions was developed by Wiener\cite{Wiener1938}, who showed that a Gaussian process with a finite variance can be expanded as an infinite series of Hermite polynomials of Gaussian random variables with deterministic series coefficients.   
This approach was popularized in the 1990s by Ghanem and co-workers as a numerical method for propagating uncertainty in stochastic finite element computations \cite{Ghanem2003}. Shortly thereafter, these methods were expanded to problems involving non-Gaussian random variables by Xiu and Karniadakis \cite{Xiu2002,Xiu2010} with the idea that for non-Gaussian random variables, a similar expansion can be performed, but the polynomial basis chosen for the expansion should be chosen such that the polynomials are orthogonal with respect to an inner product weighted by the appropriate probability density function.
By truncating such a series expansion and considering a projection onto a finite set of orthogonal polynomials, polynomial chaos methods allow a stochastic problem to be approximated by a deterministic (but higher dimensional) problem for the projection coefficients. 
Further, the orthogonality of the polynomial basis makes many numerical computations especially convenient, as the statistical moments of the output can be simplified greatly by exploiting the orthogonality, and such bases are also naturally used with Gauss quadrature methods for computing integrals which appear due to the inner products. 

In recent years, methods involving polynomial chaos expansions have also received research attention from the stochastic control community\cite{Kim2013,Mesbah2016}. Mesbah et al. \cite{Mesbah2014} developed a stochastic nonlinear model predictive control scheme with probabilistic constraints based on polynomial chaos expansions.  This was further extended by Bavdekar and Mesbah \cite{Bavdekar2016} to incorporate stochastic disturbances.  Fisher and Bhattacharya \cite{Fisher2009,Fisher2010} also developed methods of trajectory generation and solving the stochastic linear quadratic regulator problem based on polynomial chaos expansions. Boutselis et al \cite{Boutselis2019} implemented polynomial chaos expansions in a trajectory optimization method which used differential dynamic programming on the polynomial chaos expansion coefficients.  This was implemented in a receding horizon manner with probabilistic constraints by Aoyama et al \cite{Aoyama2021}. 
Similarly, Nakka and Chung \cite{Nakka2022} recently developed an iterative method based on sequential convex programming coupled with polynomial chaos expansions to solve the stochastic trajectory optimization problem with stochastic disturbances and probabilistic constraints. 

Here we introduce the preliminary material regarding generalized polynomial chaos expansions (gPC), as they apply to stochastic control systems.  That is, an ordinary differential equation with a deterministic control input, where the stochasticity is associated with uncertainties in the initial conditions and constant parameters.  The description here largely follows the accounts of Xiu \cite{Xiu2010} and Boutselis \textit{et al.}\cite{Boutselis2019}. 
Consider the stochastic control system 
\begin{equation} \label{eq:sde}
    \frac{d}{dt}x(t;Z) = f(x(t;Z),u(t),t;\theta(Z)), \qquad x(t_0;Z) = x_0(Z)
\end{equation}
where $\theta (Z)$ are model parameters, $x_0(Z)$  is the initial state, $x\in\RR^n$ is the system state, $u\in\RR^{n_c}$ is the control input,  and $Z\in\RR^d$ is a continuous random variable with a probability density function $\rho(z)$ .
Throughout this paper, we will use the shorthand notations
\[
\begin{split}
    \<g\> &= \EE[g(Z)] = \int f(z) \rho(z) dz\\
    \<g,h\> &= \EE[g(Z)h(Z)] = \int g(z)h(z)\rho(z) dz
\end{split}
\] 
and so forth.
Let $\{\phi_k(Z)\}_{k=0}^K$ be a set of gPC basis functions, that is, polynomials satisfying the following orthogonality condition with respect to a probability density $\rho(z)$
\begin{equation}
\EE[\phi_i(Z)\phi_j(Z)] = \<\phi_i,\phi_j\> = \int \phi_i(z) \phi_j(z) \rho(z) dz = \gamma_i\delta_{ij}
\end{equation}
where $\gamma_i = \<\phi_i^2(Z)\>$.  If we choose a basis of all gPC polynomials up to degree $r$, then $K = \frac{(r+d)!}{r!d!}-1$.
We then approximate the parameters, the initial condition, and the solution by their gPC expansion in polynomials of the random variable $Z$
\begin{subequations}
    \begin{align}
    \theta(Z) &\approx \sum_{j = 0}^K \hat{\theta}_j \phi_j(Z) \\
    x_i(t_0;Z) &\approx \sum_{j = 0}^K  \hat{x}_{ij}(t_0)\phi_j(Z) \\
    x_{i}(t;Z) &\approx \sum_{j = 0}^K \hat{x}_{ij}(t)\phi_j(Z)      \label{eq:gpcProjSol}  
    \end{align}
\end{subequations}
where $i = 1, \dots, n$ is a state coordinate index and the deterministic coefficients are given as follows. 
\[
   \hat{\theta}_j = \frac{\<\theta,\phi_j\>}{\<\phi_j^2\>}, \qquad 
   \hat{x}_{ij}(t_0) = \frac{\<x_i(t_0),\phi_j\>}{\<\phi_j^2\>}, \qquad 
   \hat{x}_{ij}(t) = \frac{\<x_i(t),\phi_j\>}{\<\phi_j^2\>} \qquad 
\]
It is shown by Xiu \cite{Xiu2010} that these expansions give the best approximation in $\PP^d_K$, the space of polynomials of $Z\in\RR^d$ of degree up to $K$.  However, at this point, the approximation is of little practical use, as Eq. \ref{eq:gpcProjSol} depends on the unknown solution $x(t,Z)$.
To resolve this, we apply a stochastic Galerkin method, wherein we seek an approximate solution in $\PP_K^d$ such that the residual is orthogonal to the space $\PP_K^d$ \cite{Xiu2010}. 
Substituting these approximations into the Eq. \ref{eq:sde} on the left hand side, we have
\begin{equation}
\frac{d}{dt}x_i(t;Z)\approx
\frac{d}{dt}
\left[ 
\sum_{j=0}^K \hat{x}_{ij}(t) \phi_j(Z)\right] =  
\sum_{j=0}^K
\frac{d\hat{x}_{ij}}{dt}\phi_j(Z) 
\end{equation}
and on the right hand side
\begin{equation}\label{eq:sde_rhs_gpcapprox}
    f(x(t;Z),u(t),t;\theta(Z))
    \approx
    f\left(\sum_{j=0}^K \hat{x}_{ij}(t) \phi_j(Z),u(t),t,\sum_{j = 0}^K \hat{\theta}_j \phi_j(Z)\right)\
\end{equation}
To perform the Galerkin projection, we multiply each side by $\phi_k(Z)$ with $k=1,\dots, K$ and apply the inner product on each side.  Exploiting the orthogonality of the basis functions, this yields the following deterministic ODE for the time evolution of the projection coefficients
\begin{equation} \label{eq:gpc_ode}
\frac{d\hat{x}_{ij}}{dt}
= \frac{\<f_i(x(t;z),u(t),t;\theta(z))\,,\,\phi_j(z)\>}{\<\phi_j(z)\,,\,\phi_j(z)\>} 
\end{equation}
where $f$ is also approximated as given in Eq. \ref{eq:sde_rhs_gpcapprox}.
In practice, integrals of this form will be evaluated by Gaussian quadrature. 
Denoting the vector of gPC coefficients as 
\[X = [\hat{x}_{11},\, \hat{x}_{12},\, \dots ,\, \hat{x}_{1K},\, \hat{x}_{21},\, \dots ,\, \hat{x}_{2K},\,\dots ,\, \hat{x}_{nK}]^T\] 
Eq. \ref{eq:gpc_ode} can be rewritten in the following, more compact form. 
\begin{equation}  \label{eq:gpc_ode_compact}
    \frac{dX}{dt} = \mathbf{f}(X,u)
\end{equation}

Once the time evolution of the projection coefficients $\hat{x}_{ij}$ is found, the gPC approximations of the first and second (central) moments of the states, $\mu$ and $\sigma$ can be recovered as follows: 
\begin{subequations} \label{eq:moments}
\begin{align}
\mu_i &= \int x_i(z)\rho(z) dz \approx \sum_{j=0}^K \hat{x}_{ij}\<\phi_j\> = \hat{x}_{i0}  \label{eq:gpc_mean}\\
\sigma_{ij} &= \int (x_i(z)-\mu_i)(x_j(z) - \mu_j)\rho(z) dz  = \<x_ix_j\>-\mu_i\mu_j\\
&\approx \sum_{k=0}^K\hat{x}_{ik}\hat{x}_{jk}\<\phi_k^2\>-\mu_i\mu_j = \sum_{k=1}^K\hat{x}_{ik}\hat{x}_{jk}\<\phi_k^2\>\nonumber.
\end{align}
\end{subequations}
where in the last expression we have used $\mu_i\mu_j \approx \hat{x}_{i0}\hat{x}_{j0}$ as in Eq. \ref{eq:gpc_mean}.
Higher order moments could be determined in a similar fashion.

\subsection{Optimal control formulation}  \label{sec:gpc_optimalcontrol}
We now formulate an optimal control problem to steer an initial density to a target destination using the gPC expansion to model the density transport.  
To do this, we will consider an optimal tracking problem in terms of the first and second moments, as given in Eq. \ref{eq:moments}. 
Denoting the expressions for the relevant moments as
\begin{equation}
    y = M(X)\:,
\end{equation}
we seek to optimize a cost which is quadratic in the error from the target moments, $\yr$, and the control effort, $u$, as expressed by the in-horizon and terminal cost functions, $l$ and $l_H$ given below
\begin{subequations} \label{eq:oc_loss}
\begin{align}
    l(X,u) &= (M(X)-\yr_t)^TS(M(X)-\yr) + u^TRu\\[1ex]
    l_H(X) &= (M(X)-\yr)^TS_H(M(X)-\yr)
\end{align}
\end{subequations}
where $S$, $R$, and $S_H$ are weighting matrices penalizing the in-horizon moment errors, the control effort, and the moment error at the terminal state, respectively. 
To solve this numerically, we write this as a discrete-time optimal control problem as follows
\begin{subequations}
\begin{align}
    \min_{u_1,u_2,\dots,u_{H-1}} &\sum_{t=1}^{H-1}l(X_t,u_t) + l_H(X_H)\\[1ex]
    \mathrm{s.t.} \qquad & X_{t+1} = F(X_t,u_t) \label{eq:disc_dyn2} \\
    & y_t = M(X_t)
\end{align}
\end{subequations}
where $H$ is the number of timesteps in the time horizon, the subscripts denote the discrete time step, and Eq. \ref{eq:disc_dyn2} represents the discrete time version of Eq. \ref{eq:gpc_ode_compact}.
We solve this problem numerically using differential dynamic programming  \cite{Jacobson1968, Yakowitz1984, Tassa2012}.
Differential dynamic programming (DDP) computes a locally optimal control around a nominal trajectory by minimizing a quadratic approximation of the value function along this trajectory, and then iteratively optimizing about the new trajectories obtained by applying the locally optimal control. 
A brief review of this procedure, along with the needed gradients in terms of the gPC expansions are given in \ref{sec:ddp}.

\section{Mobile rotors with velocity control} \label{sec:movingrotors}
In this section, we apply the method described in Section \ref{sec:gpc} to the problem of steering a distribution of fluid particles using the flow field produced by a group of microrotors, where the rotor translational velocities can be controlled directly.  

\begin{figure}
    \centering
    \includegraphics[width = 0.45\linewidth]{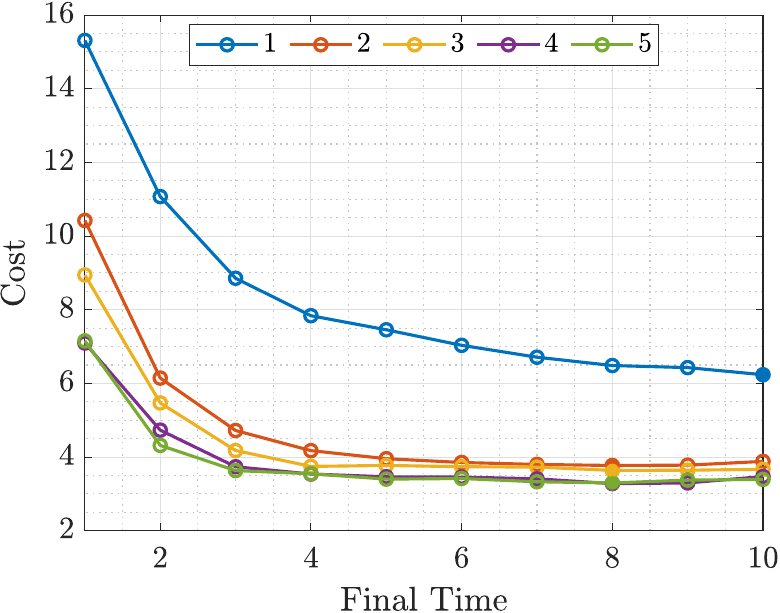}
    \caption[Comparison of the optimal control cost for steering a density to a target using mobile microrotors, as described in Sec. \ref{sec:movingrotors}, for varying final time and varying number of rotors.]{Comparison of the optimal control cost for steering a density to a target using mobile microrotors, as described in Sec. \ref{sec:movingrotors}, for varying final time (horizontal axis) and varying number of rotors (curves). The minimum point on each curve is shown by a filled circle. }
    \label{fig:costcomp}
\end{figure}

With this system we consider the problem of steering an initial particle distribution described by the Gaussian density, $\N([1,1],0.025I_2)$ to a target mean set to $(-1,-1)$ and the target variance zero.  That is, we seek to minimize the variance while steering the distribution toward the target. 
The initial positions of the rotors are set to be evenly distributed along a circle of radius $0.2$ about the target mean position.  This forces the rotors to move away from the target, retrieve the distribution of particles, and steer it to the target.  
In the numerical solution to the optimal control problem, we use a timestep of $\Delta t = 0.01$ and a cost function of the form in Eq. \ref{eq:oc_loss}, 
where the moments of interest are $M(X_t) = [\mu_1,\mu_2,\sigma_{11},\sigma_{22}]^T$ and the cost function weight matrices are chosen as follows.  
\begin{subequations}
\begin{align}
    S &= \mathrm{diag}([0.1,0.1,0.1,0.1])\Delta t\\[1ex]
    S_H &= \mathrm{diag}([1000,1000,1000,1000])\Delta t\\[1ex]
    R &= \mathrm{diag}([\underbrace{1,\dots, 1}_{n_r},\underbrace{0.1,\dots, 0.1}_{n_r},\underbrace{0.1,\dots, 0.1}_{n_r}])\Delta t
\end{align}
\end{subequations}
These cost weights are chosen to place a heavy penalty on the error at the terminal state, while penalizing the error in-horizon relatively little.  This allows the optimizer the freedom to apply controls within the horizon that could result in higher cost, so long as the terminal cost is minimized.  The cost on the rotor strengths is set to $1$, giving a scale for the terms in the cost function and the penalty on the magnitude of the translational velocity is set to $0.1$. 
We vary the number of rotors as $n_r = 1,2,3,4,5$ and $6$ and consider final times ranging from $1$ to $10$. 
The only stochastic quantities considered are the two random variables describing the initial particle distribution.  Since this initial distribution is a two dimensional Gaussian, the gPC basis used consists of Hermite polynomials of $x_p$ and $y_p$.  In all of the results shown here, we consider gPC expansions of polynomials of degree up to $3$.   
This optimal control problem is solved numerically using the differential dynamic programming scheme detailed in \ref{sec:ddp}.
A comparison of the optimal control cost for varying final time and varying numbers of rotors is shown in Fig. \ref{fig:costcomp}.  The costs here are evaluated on the true moments, as obtained by applying the control obtained from DDP on the gPC expansion in a Monte Carlo simulation of $10^4$ particles sampled from the initial distribution.

Fig. \ref{fig:costcomp} shows that there is a significant improvement in the performance of this controller for two rotors as compared to one rotor, but the improvement diminishes incrementally as more rotors are added. For numbers of rotors greater than $n_r=4$, we see negligible improvement in the cost. In comparing the performance as a function of the final time, in each curve the largest costs occur for the shortest time horizons considered, with the cost decreasing gradually with increasing final time.  In each of the cases with $n_r\geq2$, the curve reaches a minimum point at $t_f=8$, whereas for $n_r=1$, the lowest value obtained from this set of simulations occurs at the largest final time considered, $t_f = 10$.  The large costs at small $t_f$ values is mostly due to the large error in the terminal state, as the rotors are not able to drive the distribution to the target in such a short time without expending a large control effort or without stretching the distribution very significantly.  As the final time increases, eventually a minimizing solution is found for each number of rotors, such that the error in the moments does not improve even if more time is allowed.  Beyond this point, the more cost is accumulated within the horizon without any improvement of the terminal cost, and therefore the total cost begins to increase very slowly for higher final times, as is seen in Fig. \ref{fig:costcomp} for $n_r\geq 2$ and $t_f>8$. 
For these reasons, in the remainder of this section, we will focus in detail on the case of four rotors, $n_r = 4$ and $t_f = 8$, as this is in some sense, a best case.  For larger numbers of rotors the cost improvement is seen to be negligible and for longer time horizons the cost is seen to increase. 

\begin{figure*}
  \hspace{1.3cm}$n_r=1$  \hspace{1.8cm} $n_r=2$ \hspace{1.8cm} $n_r=3$ \hspace{1.8cm} $n_r=4$ \hspace{1.8cm} $n_r=5$ 
  \vspace{-1ex}
\begin{center}
    \includegraphics[height = 3.34cm,trim={45 0 90 17},clip]{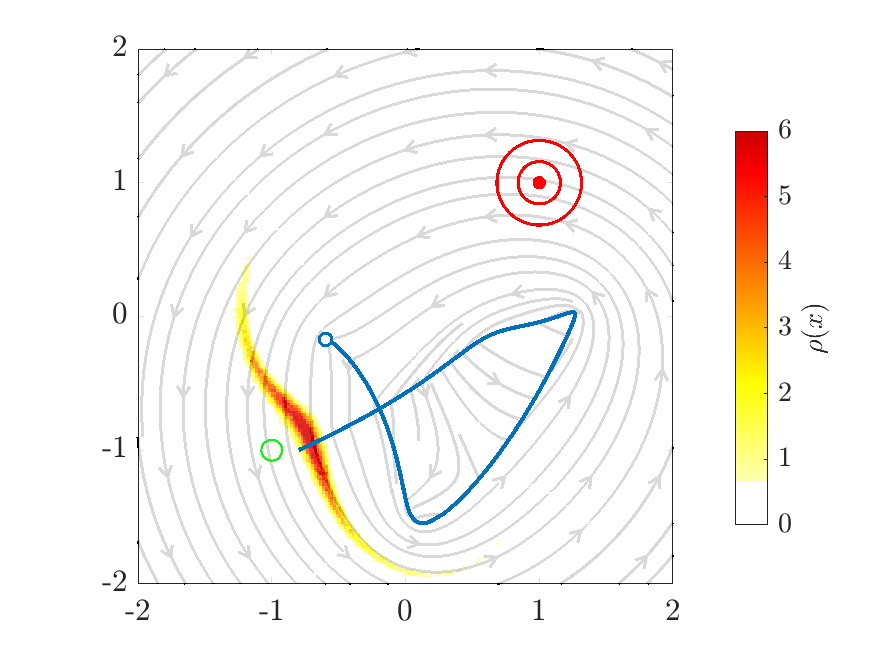}
    \includegraphics[height = 3.34cm,trim={65 0 90 17},clip]{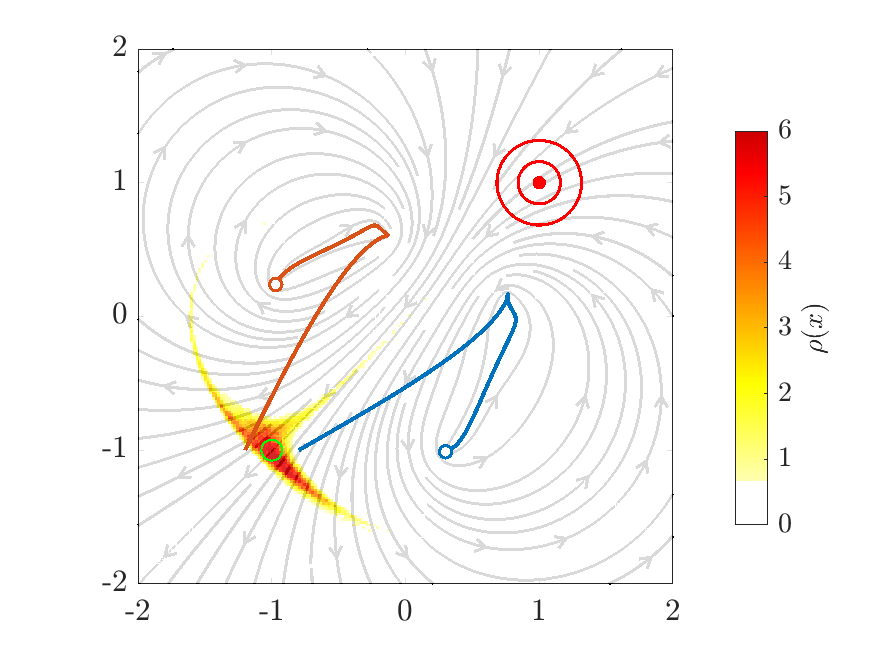}
    \includegraphics[height = 3.34cm,trim={65 0 90 17},clip]{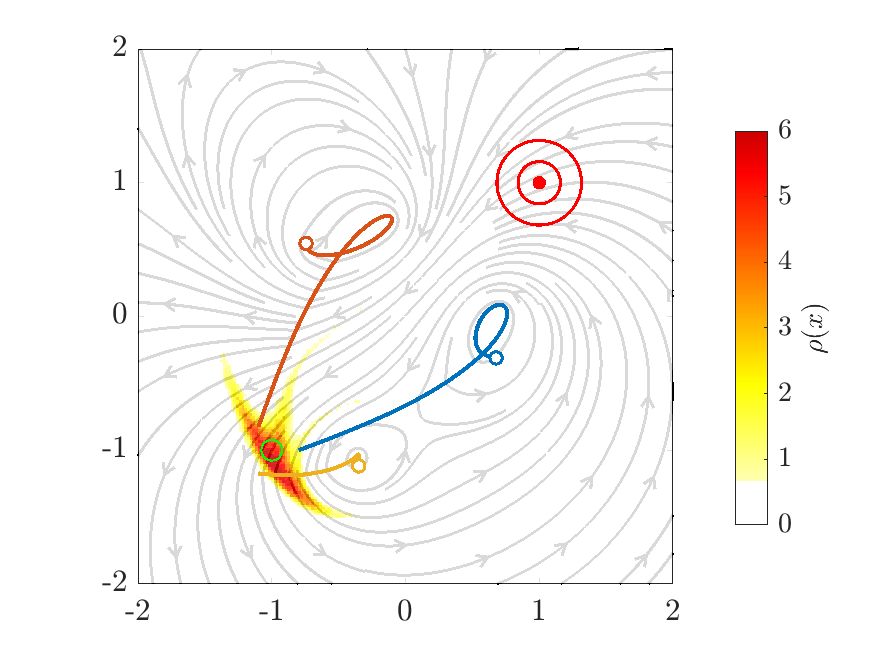}
    \includegraphics[height = 3.34cm,trim={65 0 90 17},clip]{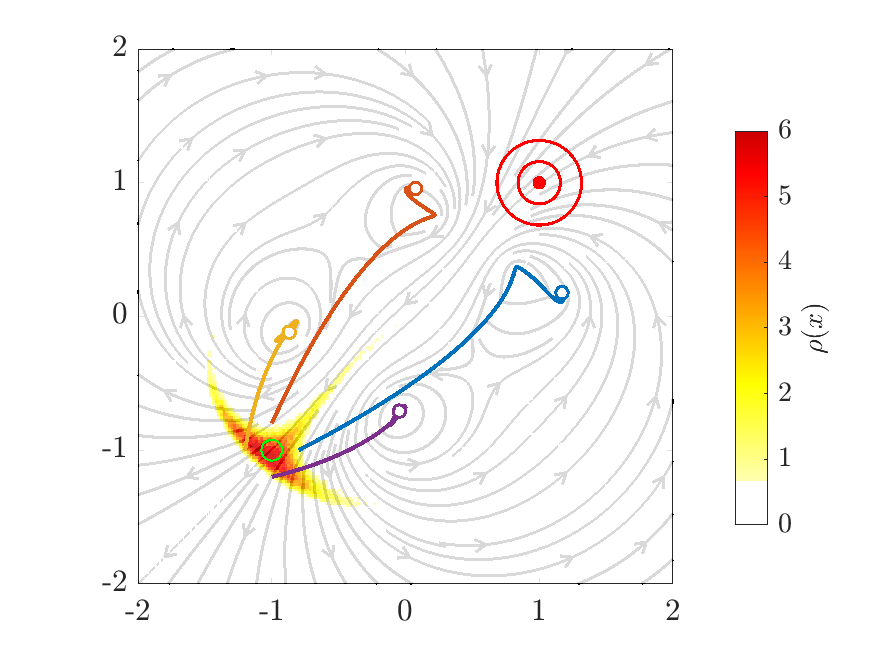}
    \includegraphics[height = 3.34cm,trim={65 0 20 17},clip]{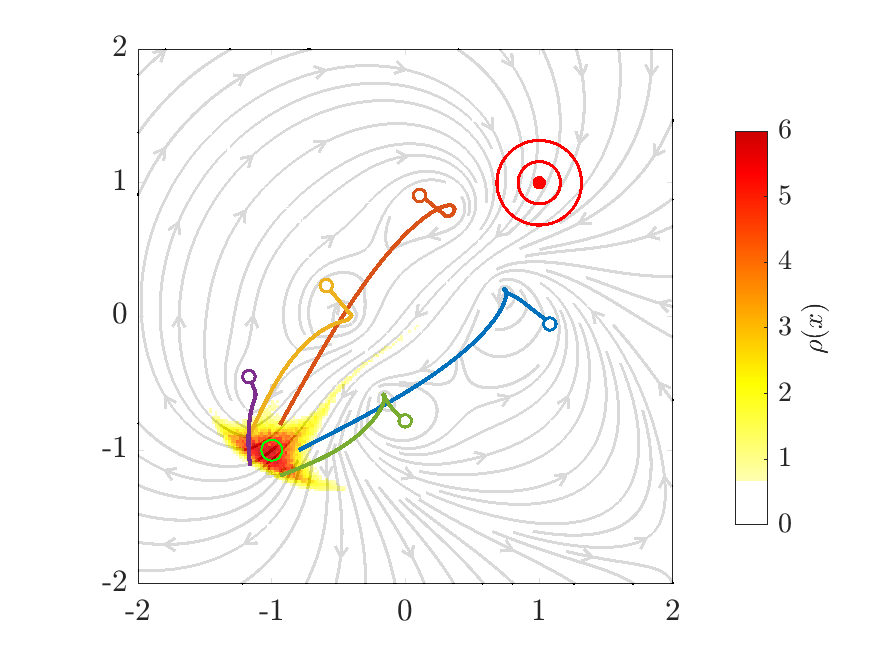}
    \end{center}
    \caption{Streamlines show the time-averaged velocity field for the optimal case for varying number of rotors.  Color depicts the particle density at the final time for each case.  The red marker and contours depict the mean, 1$\sigma$ and 2$\sigma$ contours for the initial fluid particle distribution.  Colored lines indicate the paths of the rotors. The green circle indicates the target at (-1,-1).}
    \label{fig:movrot_vel_tavg}
\end{figure*}

Figure \ref{fig:movrot_vel_tavg} depicts the solution corresponding to the lowest cost point for each number of rotors shown in Fig. \ref{fig:costcomp}. The streamlines shown represent the time-averaged velocity field, averaged over the full simulation.  The color in each figure shows the probability density of the particle distribution at the final time for each case.  The colored lines indicate the paths of each rotor in the simulation, and the red marker and lines show the mean, $1\sigma$, and $2\sigma$ contours of the initial particle density. In each case, we see that the rotors move away from their starting points at the target, in order to drive the the particle distribution toward the goal with less rotor strength.  In all of the cases with multiple rotors, the rotors split roughly symmetrically about the line between the target and the initial distribution location. The rotors on either side of this line then apply torques of nearly the same magnitude, but of opposite direction, which pulls the distribution directly along the line separating the rotors. This form of the solution is seen very clearly in Fig. \ref{fig:movrot_vel_tavg} for the cases with even numbers of rotors, $n_r = 2$ and $n_r = 4$.  The final densities shown in these plots also indicate that there is considerably less stretching of the particle distribution for increasing number of rotors. 

\begin{figure*}
    \hspace{1.3cm}$t = 0.0$  \hspace{1.8cm} $t = 2.0$ \hspace{1.8cm} $t = 4.0$ \hspace{1.8cm} $t = 6.0$ \hspace{1.8cm} $t = 8.0$ 
    \begin{center}
    \includegraphics[height = 3.34cm,trim={45 0 90 17},clip]{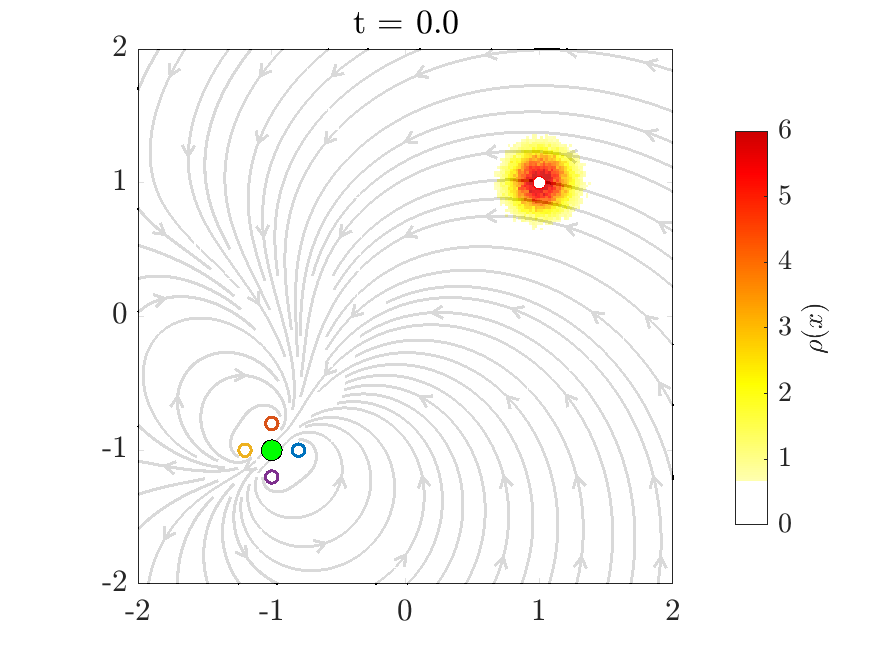}
    \includegraphics[height = 3.34cm,trim={65 0 90 17},clip]{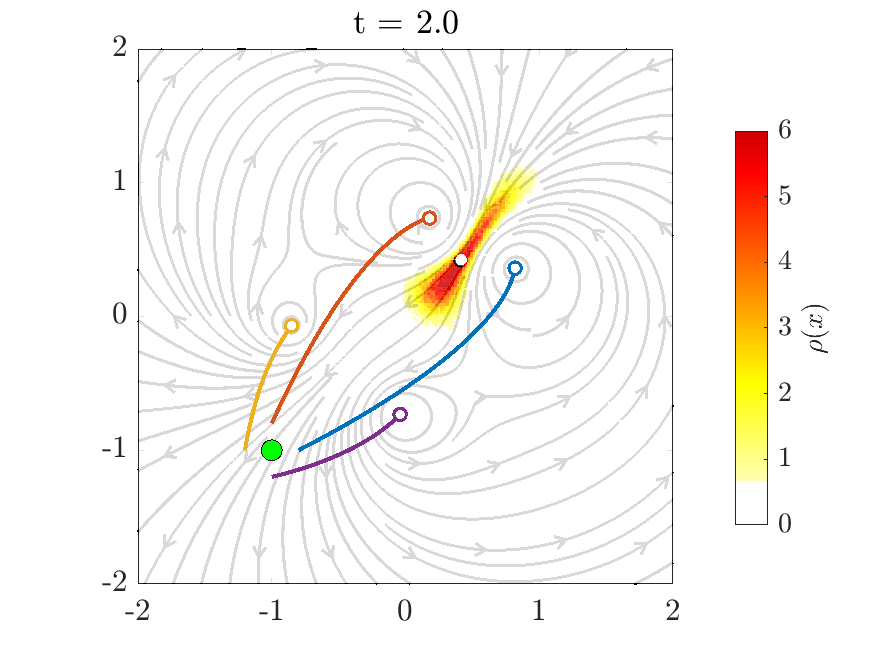}
    \includegraphics[height = 3.34cm,trim={65 0 90 17},clip]{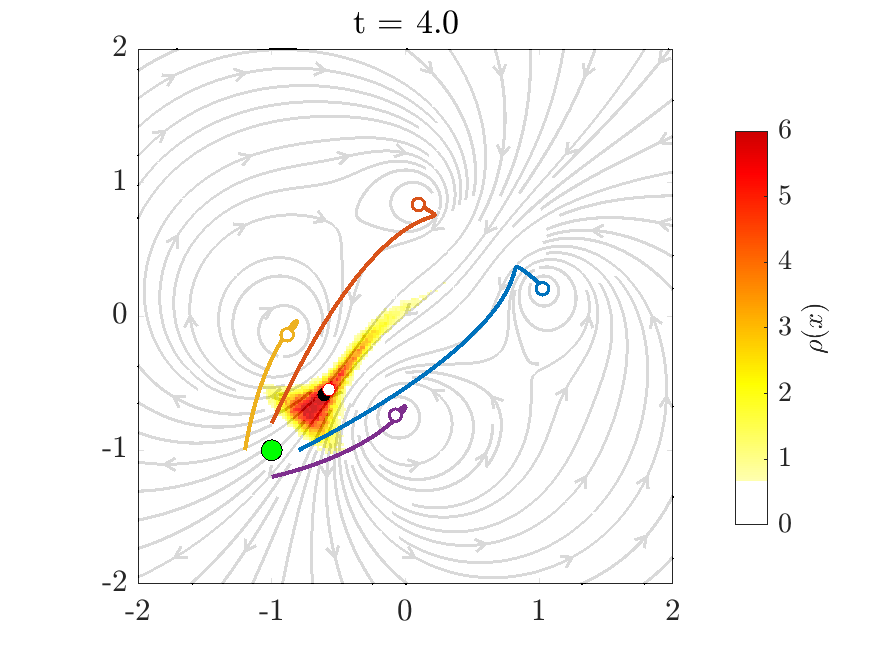}
    \includegraphics[height = 3.34cm,trim={65 0 90 17},clip]{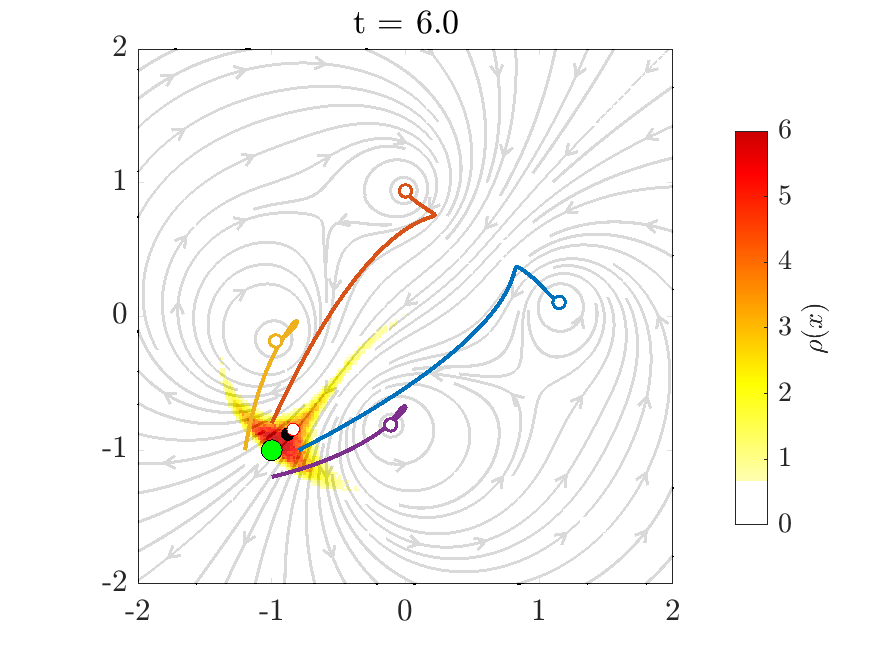}
    \includegraphics[height = 3.34cm,trim={65 0 20 17},clip]{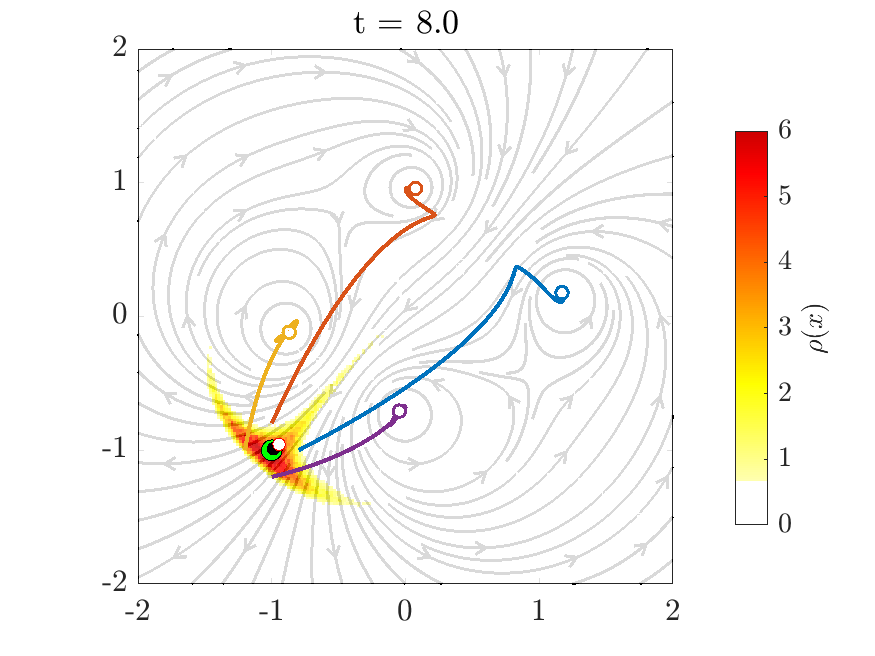}
    \caption[Optimal solution for the case of four rotors, $n_r = 4$ and final time $t_f = 8$ shown as snapshots from the time sequence.]{Optimal solution for the case of four rotors, $n_r = 4$ and final time $t_f = 8$ shown as snapshots from the time sequence.  Streamlines show the direction of the instantaneous velocity fields at each instant shown.  Color depicts the fluid particle density at each time instant.  Colored lines show the rotor paths, with circles indicating the rotor position at each instant. Black circle shows the sample mean, white circle shows the gPC-predicted mean. }\label{fig:nr4_traj}
    \end{center}
    \vspace{-2em}
\end{figure*}

\begin{figure}
    \centering
    \includegraphics[width=0.45\linewidth]{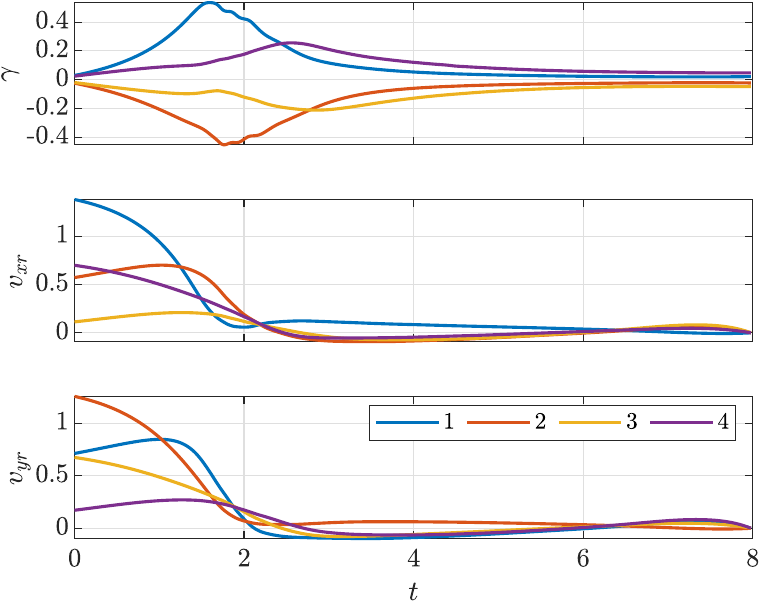}
    \vspace{-1ex}
    \caption[Controls computed from the DDP procedure for the case of $n_r=4$ and $t_f=8$.]{Controls as computed from the DDP solution for the case of $n_r=4$ and $t_f=8$.  Each curve represents a different rotor, with colors corresponding to the colored trajectories in Fig. \ref{fig:nr4_traj}. Top: rotor strengths, $\gamma$. Middle:  $x$-component of rotor translational velocity. Bottom: $y$-component of rotor translational velocity.  Rotor numbering goes counterclockwise from the rotor directly right of the target in the initial configuration. }
    \label{fig:nr4_control}
\end{figure}

\begin{figure}
    \centering
    \includegraphics[width=0.35\linewidth]{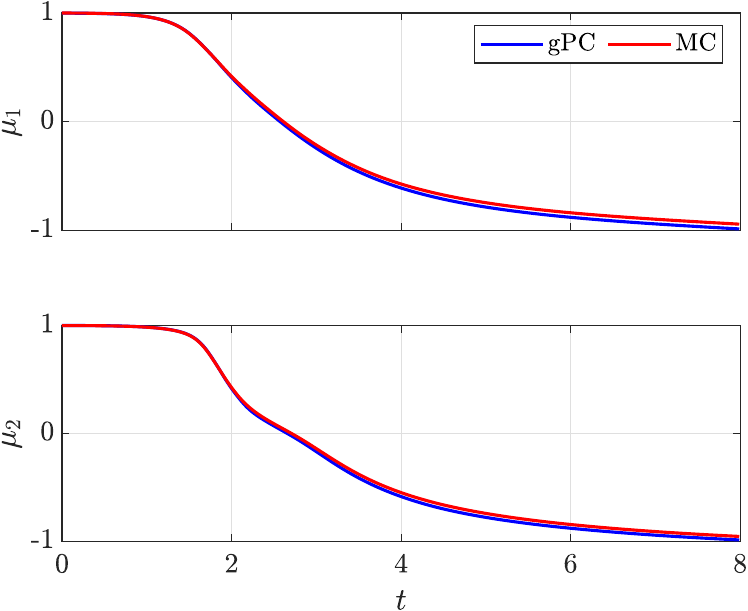}
    \includegraphics[width=0.35\linewidth]{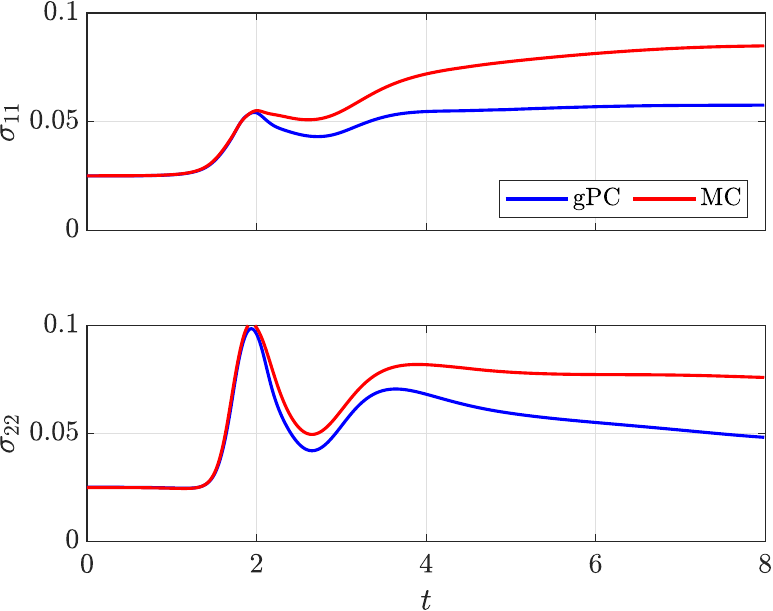}
    \caption[Moment propagation by the polynomial chaos expansion method (labelled gPC) as compared to the Monte Carlo simulation for the optimal solution of the case of $n_r=4$.]{Moment propagation by the polynomial chaos expansion method (labelled gPC) as compared to the Monte Carlo simulation of $10^4$ points sampled from the initial particle distribution (labelled MC) for the optimal solution of the case of $n_r=4$, as shown in Fig. \ref{fig:nr4_traj}.  Left column: components of the mean, $\mu$. Right column: components of the covariance, $\sigma$.    }
    \label{fig:nr4_moments}
\end{figure}

The solution for the case of four rotors, $n_r = 4$, and final time, $t_f = 8$ is shown in more detail in Figs. \ref{fig:nr4_traj}, \ref{fig:nr4_control}, and \ref{fig:nr4_moments}. Fig. \ref{fig:nr4_traj} shows a sequence of snapshots from this solution, showing the particle distribution, rotor location, and rotor trajectories at five instances between $t=0$ and $t=8$.  Fig. \ref{fig:nr4_control} shows the rotor strengths and translational velocity controls used to produce this solution.  The plot of the strengths in the top row of Fig. \ref{fig:nr4_control} also highlights the symmetry of the strengths for rotors on opposite sides of the line connecting the initial distribution location to the target point mentioned in the previous paragraph, as rotors $2$ and $3$ take nearly the negative strengths of rotors $1$ and $2$.  
From these plots, it is clear that the majority of the transport occurs early in the solution, but after the rotors have had enough time to move toward the distribution, approximately from $t = 1$ to $t=3$.  This can also be seen in Fig. \ref{fig:nr4_moments}, which shows the moment propagation results for this case. This figure compares the moment propagation using the polynomial chaos expansion in polynomials of up to degree $3$ with the Monte Carlo simulation of $10^4$ points sampled from the initial distribution. From this, we see close agreement of the moments, especially in the first half of the trajectory.  

\clearpage

\clearpage

\section{Analysis of Lagrangian transport} \label{sec:lcs} 
We will now further analyze the unsteady flow field associated with the optimal control solution by applying tools from the study of nonautonomous dynamical systems. Specifically, we will compute the Lagrangian coherent structures (LCS) associated with this flow field by analyzing the ridges of the finite-time Lyapunov exponent (FTLE) field.  The FTLE/LCS method is a geometric method for the study of transport properties of a dynamical system, specifically developed as an extension of phase space methods for autonomous systems to nonautonomous systems, as the flow structures identified as ridges of the FTLE field can be thought of as analogous to separatrices or  stable and unstable manifolds in an autonomous system. 
Previous works have studied the relationship between the optimal control of particle motion in an unsteady flow and the Lagrangian coherent structures of the flow \cite{Inanc2005,Senatore2008,Ramos2018,Krishna2022, Krishna2023}.  These works showed that steering a particle (or underwater vehicle in those cases) to follow the coherent structures of the background flow which tend to guide the particle to the target destination leads to energy-efficient navigation. Krishna \textit{et al.} \cite{Krishna2023} studied how coherent structures are modified by the addition of an optimal feedback control policy to a flow field. 
Rather than studying the control of individual agents or particles in an uncontrolled flow as in those works, here we have developed an optimal control strategy to control a group of microrotors to generate a flow which transports ensembles of fluid particles to a target. 
By studying the LCS produced by this optimized flow field, we see how the optimal control generates flow structures which entrain the fluid particles of interest and direct them to the target.  

Below we briefly review the derivation and uses of the finite time Lyapunov exponent, largely following Shadden \cite{Shadden2005,Shadden2005a}.  See Refs. \cite{Haller2015, Allshouse2015} for more recent reviews of such methods. 

The FTLE can be derived by quantifying the stretching induced by the flow map, $\Phi_{t_0}^{t_0+\tau}$ of a time-varying system for a finite timespan $\tau$ from an initial time $t_0$. If the initial displacement between two points is $\delta x(t_0)$, then after time $\tau$, the perturbation grows (or decays) to 
\[
\begin{split}
\delta x(t_0+\tau)
&= 
\Phi_{t_0}^{t_0+\tau}(x(t_0) + \delta x(t_0)) - \Phi_{t_0}^{t_0+\tau}(x(t_0))
\\
&= \frac{d\Phi_{t_0}^{t_0+\tau}}{dx}\bigg|_{x=x(t_0)}\delta x(t_0) + \mathcal{O}(\|\delta x(t_0)\|^2)\,.
\end{split}
\]
So, neglecting the higher order terms, the leading order of the magnitude of the perturbation after time $\tau$ is 
\begin{subequations}
    \begin{align}
           \|\delta x(t_0+\tau) \| 
    &= 
    \sqrt{\left(\frac{d\Phi_{t_0}^{t_0+\tau}}{dx}\bigg|_{x=x(t_0)}\delta x(t_0)\right)^{T}\left(\frac{d\Phi_{t_0}^{t_0+\tau}}{dx}\bigg|_{x=x(t_0)}\delta x(t_0)\right)}
    \\[1ex]
    &=
    \sqrt{
    (\delta x^{\top}(t_0)) 
    \left[
    \left(\frac{d\Phi_{t_0}^{t_0+\tau}}{dx}\right)^T
    \left(\frac{d\Phi_{t_0}^{t_0+\tau}}{dx}\right)
    \right]
    \delta x(t_0)  
    } 
    \end{align}
\end{subequations}
where the matrix in square brackets is a finite time version of the Cauchy-Green strain tensor for the flow field, denoted by $C$, 
\begin{equation}
    C = \left(\frac{d\Phi_{t_0}^{t_0+\tau}}{dx}\right)^T
    \left(\frac{d\Phi_{t_0}^{t_0+\tau}}{dx}\right).
\end{equation}
The direction of maximal stretching is given by the leading eigenvector of this matrix and the growth of a perturbation in this direction is 
\begin{equation}
    \sqrt{
    (\delta x(t_0))^{T} 
    \lambda_{\mathrm{max}}(C)
    \delta x(t_0) 
    }
\end{equation}
where $\lambda_{\mathrm{max}}(C)$ is the maximum eigenvalue of $C$. 
Then the finite time Lyapunov exponent, $\sigma_{t_0}^{\tau}$ is defined as 
\begin{equation}\label{eq:ftle}
    \sigma_{t_0}^{\tau}(x) = \frac{1}{|\tau|}\log \sqrt{\lambda_{\mathrm{max}}(C)}
\end{equation} 
so that the magnitude of the growth of a perturbation along the direction of maximum deformation can be rewritten as follows. 
\begin{equation} 
\big\|\delta x(t_0+\tau) \big\|  = \big\|\delta x(t_0)\big\|\exp{\bigg(\sigma_{t_0}^{\tau}(x)\,|\tau|\bigg)}
\end{equation}
With the definition of the FTLE in Eq. \ref{eq:ftle}, we have a definition of a scalar field which quantifies local stretching induced by the time-varying flow field at each point in space.  The ridges of high FTLE values effectively act as separatrices or transport barriers, which divide the phase space into coherent regions and through which minimal transport is allowed.  Such maximal FTLE ridges are commonly referred to as Lagrangian coherent structures. 
Just as the application of Eq. \ref{eq:ftle} with positive $\tau$ defines a scalar field which quantifies stretching in forward time, and thus ridges from which tracer particles tend to be repelled in forward time; taking a negative time increment $\tau$ quantifies stretching in backward time, and thus the associated ridges are coherent structures to which tracer particles tend to be attracted in forward time.  This effectively defines attracting and repelling Lagrangian coherent structures. 

\begin{figure} 
\centering
\begin{minipage}{0.45\linewidth}
    \centering
    \includegraphics[height=\linewidth,trim={50 0 88 0},clip]{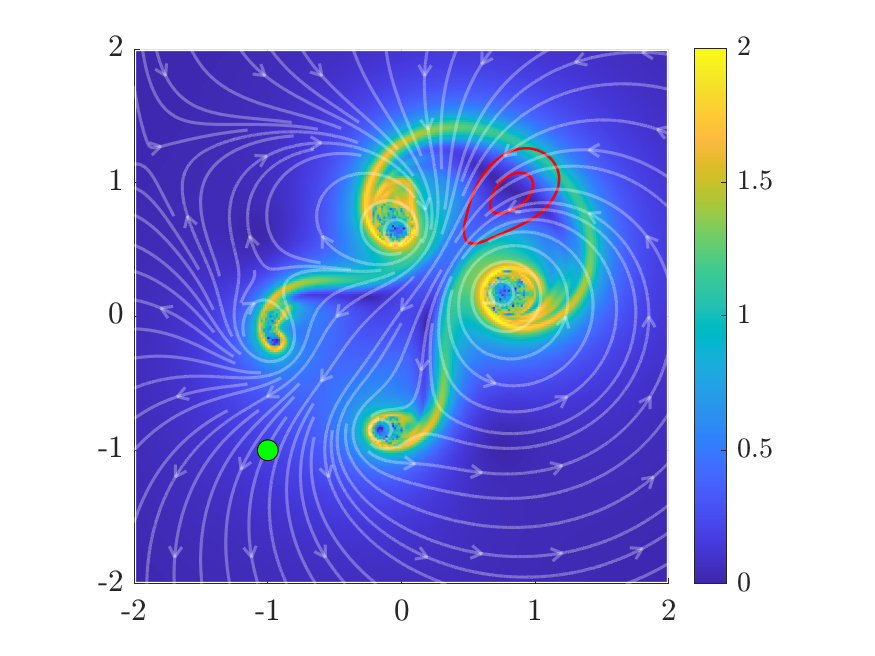}\\
    (a)
    \end{minipage}
    \begin{minipage}{0.45\linewidth}
    \centering
    \includegraphics[height=\linewidth,trim={59 0 50 0},clip]{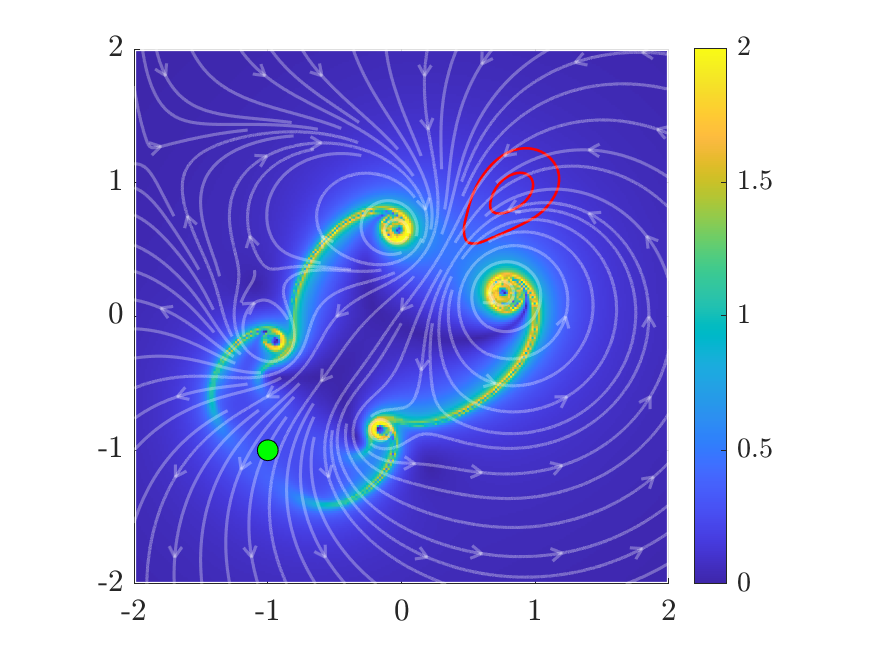} \\
    \hspace{-2em}(b)
    \end{minipage}
    \caption[FTLE field for the flow field corresponding to the optimal control for the case $n_r=4$, $t_f = 8$.]{FTLE field for the flow field corresponding to the optimal control for the case $n_r=4$, $t_f = 8$ computed at time $t = 1.5$ with FTLE integration times of $\tau = 1.5$ and $\tau = -1.5$ in panels (a) and (b) respectively.  Streamlines depict the instantaneous velocity field at $t = 1.5$.  Red contour lines show level sets of the fluid particle density, $\rho(x,1.5)$.  The level sets are shown at the values $\rho(\mu_0+\sigma_0,0)$ and $\rho(\mu_0+2\sigma_0,0)$, where $\mu_0$ and $\sigma_0$ are the mean and standard deviation of the initial density, respectively. } 
    \label{fig:ftle_movrot}
\end{figure}

With these definitions in hand, we now return to the flow field produced by the optimal control of the mobile rotor system.  As discussed previously, for the case shown in Figs. \ref{fig:nr4_traj} , \ref{fig:nr4_control}, and \ref{fig:nr4_moments}, the majority of the density transport occurs between times of $t = 0$ and $t = 3$.  This motivates studying the transport properties of the flow by computing the FTLE field for the time varying flow over this time interval. Specifically, we compute the forward and backward time FTLE field for this flow from time $t = 1.5$, with integration timespans of $\tau = 1.5$ and $\tau = -1.5$, respectively, to analyze the attracting and repelling flow structures induced by the optimal control. Computationally, this is done by integrating Lagrangian tracer particles initially distributed on a uniform grid of $250\times 250$ particles over $[-2,2]^2$ and computing the flow map Jacobians needed for the Cauchy-Green deformation tensor by central differencing.  The resulting FTLE fields are shown in Fig. \ref{fig:ftle_movrot}, along with the instantaneous streamlines and level sets of the particle density at $t = 1.5$.  The FTLE ridges here clearly reveal a structure of the time-varying flow field that is not obvious from viewing a sequence of instantaneous streamlines for the flow, as shown in Fig. \ref{fig:nr4_traj}.  Fig. \ref{fig:ftle_movrot} (a) shows that a repelling LCS is formed just behind the region of high particle concentration, as indicated by the contours of the density function.  The structure identifies the region of the fluid domain which is entrained by the rotors and pulled inward toward the target, separating this region from the outer region which remains relatively stationary in comparison. Relatedly, Fig. \ref{fig:ftle_movrot} (b) reveals attracting flow structures which effectively pull the distribution inward toward the target.


\section{Mobile rotors with torque control} \label{sec:torquecontrol}

In this section, we consider the case where the rotor velocities are not directly controlled, but instead the rotors are advected by the velocity field produced by the other rotors, as given by Eq. \ref{eq:lonly_rotorrotor}. 
With this system, we consider a similar optimal control problem to the one considered in Sec. \ref{sec:movingrotors}, where we seek to steer an initial fluid particle distribution to a desired target by controlling the rotor motion and flow fields using the rotor strengths.  
As before, we consider an initial particle distribution described by the probability density, $\N([1,1],0.025I_2)$ and a target mean set to $(-1,-1)$ and the target variance zero.  
The initial positions of the rotors are again set to be evenly distributed along a circle of radius $0.2$ about the target mean position and we use a timestep of $\Delta t = 0.01$. 
The cost function is chosen to be of the same form as before (see Eq. \ref{eq:oc_loss})
where the moments of interest are again $M(X_t) = [\mu_1,\mu_2,\sigma_{11},\sigma_{22}]^T$ and the cost function weight matrices are chosen as follows
\begin{subequations}
\begin{align}
    S &= \mathrm{diag}([0.1,0.1,0.1,0.1])\Delta t\\[1ex]
    S_H &= \mathrm{diag}([500,500,500,500])\Delta t\\[1ex]
    R &= \mathrm{diag}([\underbrace{1,\dots, 1}_{n_r}])\alpha^2\Delta t
\end{align}
\end{subequations}
where we have introduced a scaling factor $\alpha<1$.  This factor is to account for the fact that a larger magnitude of rotor strength is needed to successfully solve the control problem, as the rotor strength is responsible for defining the fluid velocity field as well as the rotor motion.
This parameter is chosen somewhat heuristically, but in practice, we find that a value of $\alpha = 1/3$ allows the optimal control problem to be solved well, and this value is used in the results shown throughout this section. 

With this formulation, we solve the optimal control problem numerically using the differential dynamic programming scheme as detailed in Sec. \ref{sec:gpc_optimalcontrol}.
Since this initial distribution is still a two dimensional Gaussian, the gPC basis used consists of Hermite polynomials of the particle position coordinates, $x_p$ and $y_p$.  In all of the results shown here, we consider gPC expansions of polynomials of degree up to $3$. 
We vary the number of rotors as $n_r = 2,3,4,5$ and consider a final time of $t_f = 10$.  We neglect the case of $n_r = 1$ here, as in this case, the rotor would be fixed, as at least two rotors must be present to generate translational motion of the rotors.  

\begin{figure}
    \centering
    \begin{minipage}{0.4\linewidth}
    \centering
    \includegraphics[width = \linewidth]{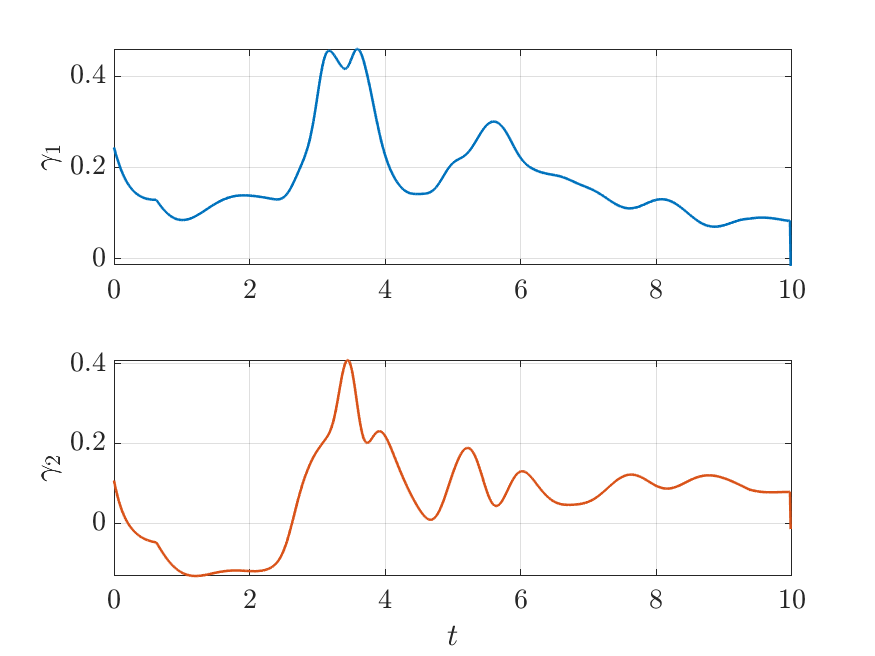}\\
    (a)
    \end{minipage}
    \begin{minipage}{0.4\linewidth}
    \centering
    \includegraphics[width = \linewidth]{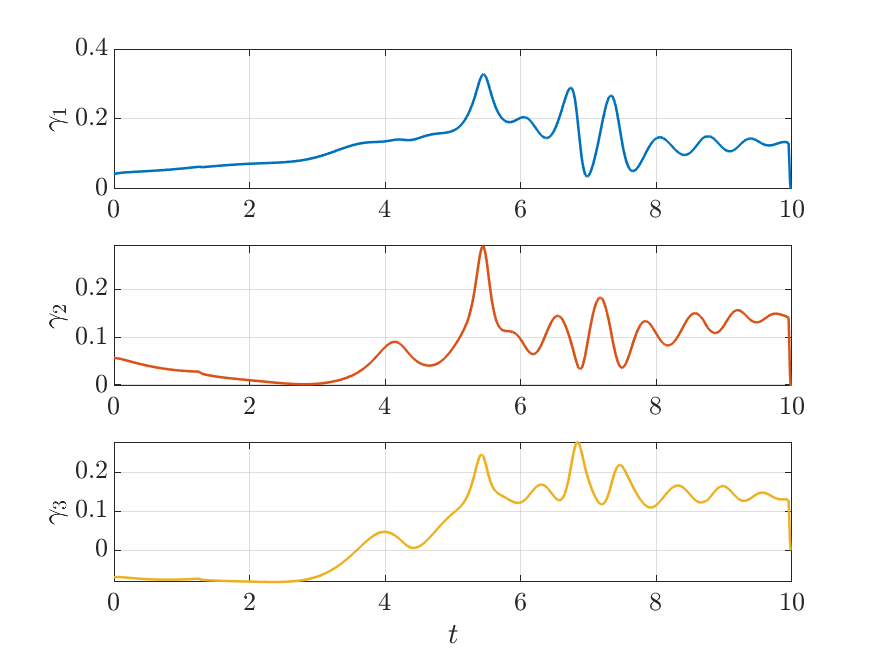}\\
    (b)
\end{minipage}
    \\
    \begin{minipage}{0.4\linewidth}
    \centering
    \includegraphics[width = \linewidth]{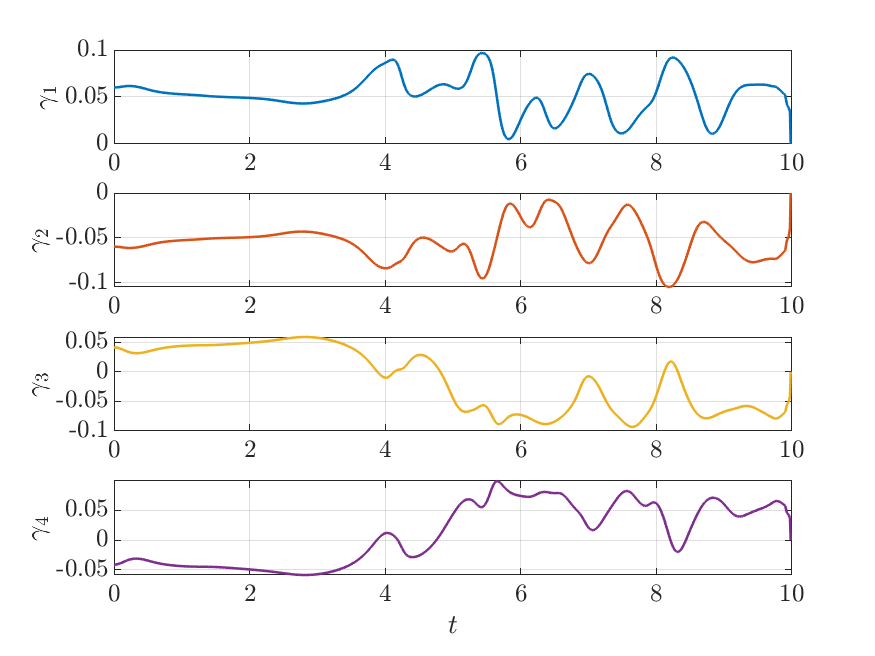}\\
    (c)
    \end{minipage}
    \begin{minipage}{0.4\linewidth}
    \centering
    \includegraphics[width = \linewidth]{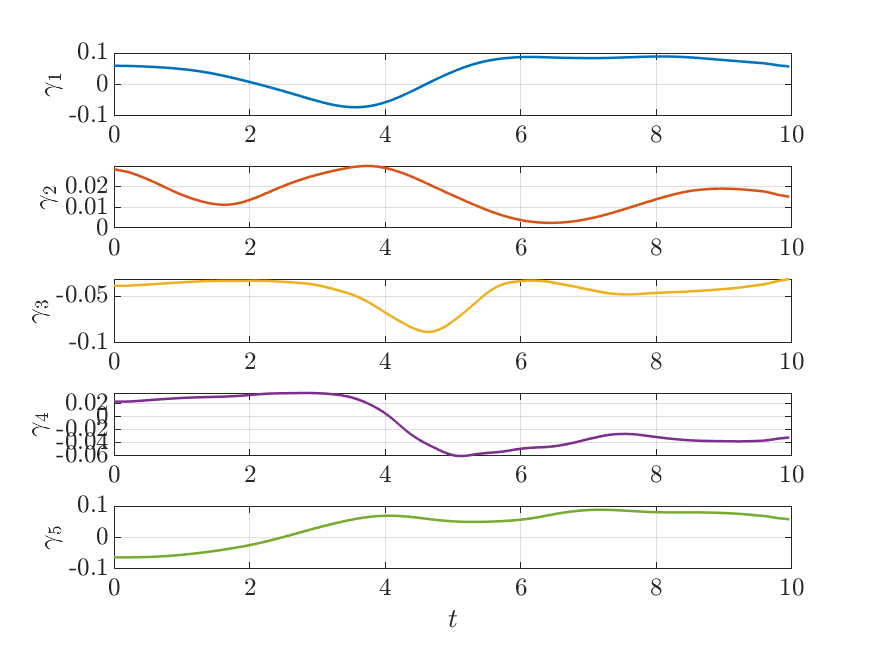}\\
    (d)
    \end{minipage}
    \\
    \caption[Control trajectories of the optimal solution for torque-only control with varying number of rotors.]{Control trajectories of the optimal solution for torque-only control with, (a) $n_r = 2$, (b) $n_r = 3$, (c) $n_r = 4$, and (d) $n_r = 5$,  corresponding to the snapshots shown in Fig. \ref{fig:lonly_traj}.  Line colors correspond to the trajectory colors in Fig. \ref{fig:lonly_traj}}.
    \label{fig:lonly_control}
\end{figure}

\begin{figure*}
\hspace{1.3cm}$t = 0.0$  \hspace{1.8cm} $t = 2.5$ \hspace{1.8cm} $t = 5.0$ \hspace{1.8cm} $t = 7.5$ \hspace{1.75cm} $t = 10.0$ 
\begin{center}
    \includegraphics[height = 3.0cm,trim={45 30 90 20},clip]{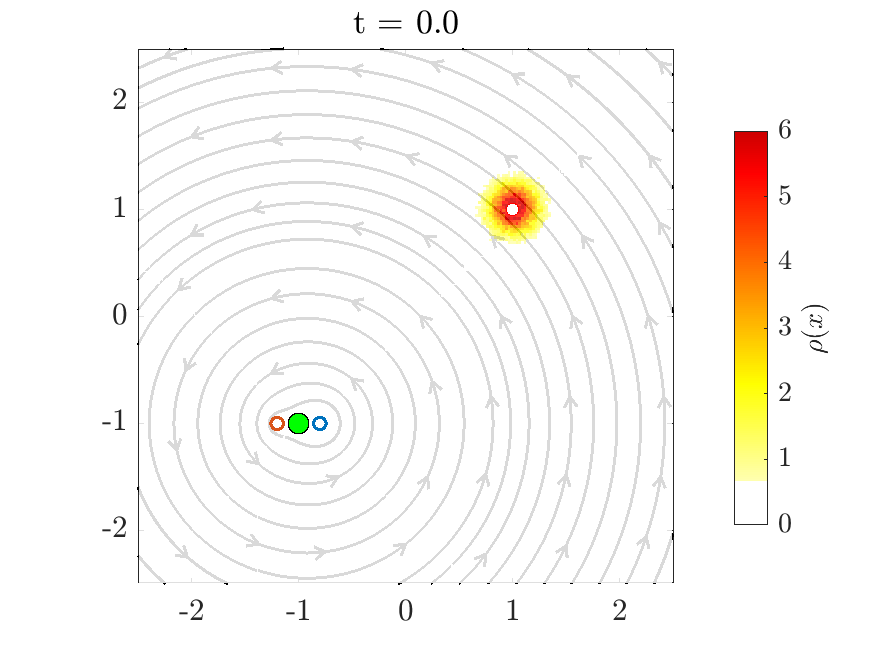}
    \includegraphics[height = 3.0cm,trim={65 30 90 20},clip]{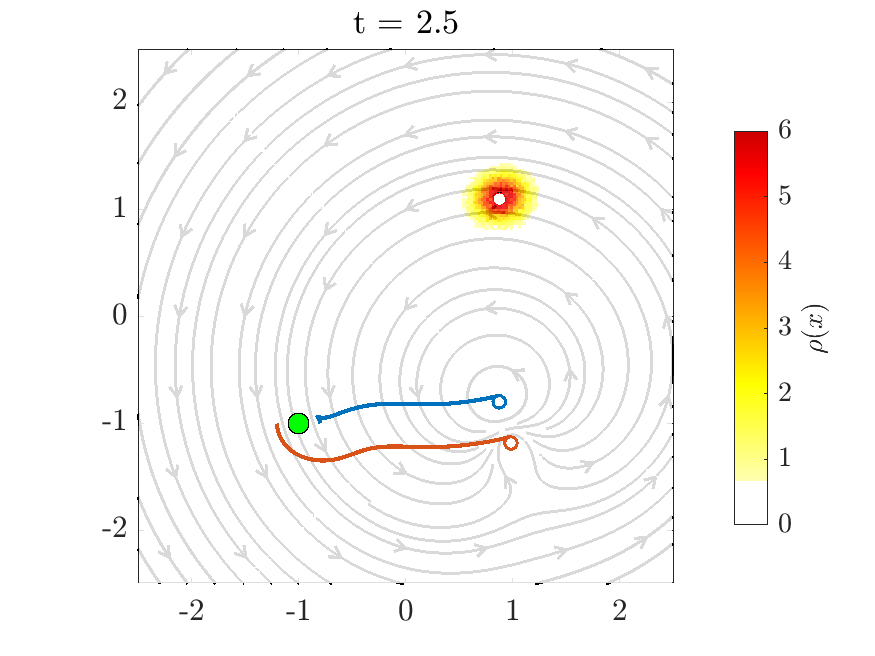}
    \includegraphics[height = 3.0cm,trim={65 30 90 20},clip]{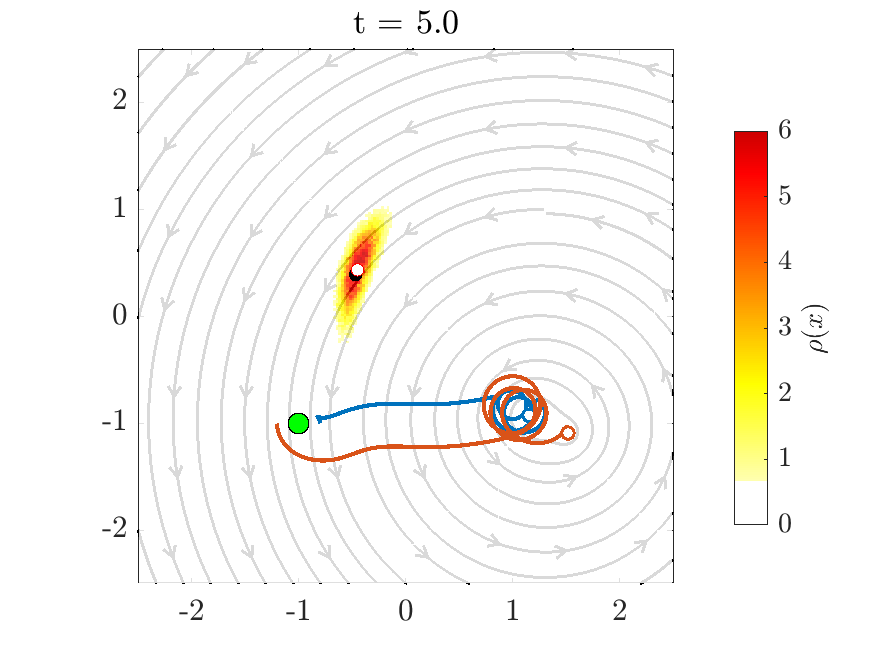}
    \includegraphics[height = 3.0cm,trim={65 30 90 20},clip]{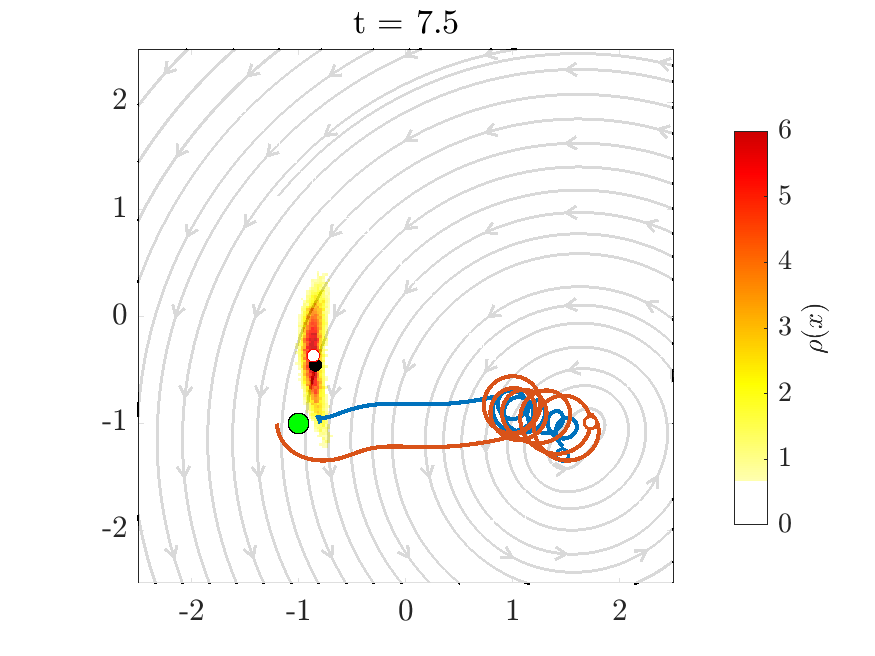}
    \includegraphics[height = 3.0cm,trim={65 30 20 20},clip]{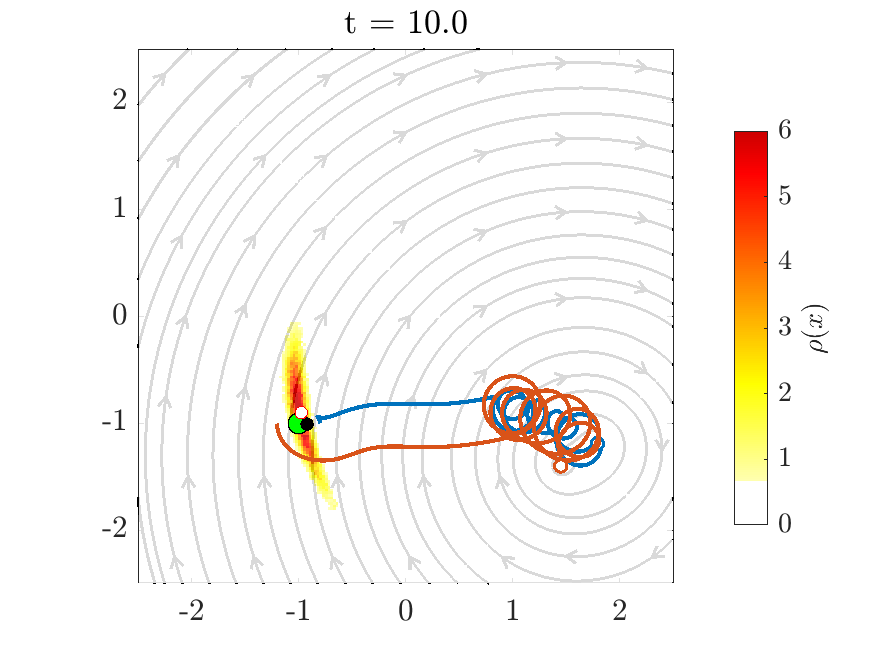}
    \\[3pt]
    \includegraphics[height = 3.0cm,trim={45 30 90 20},clip]{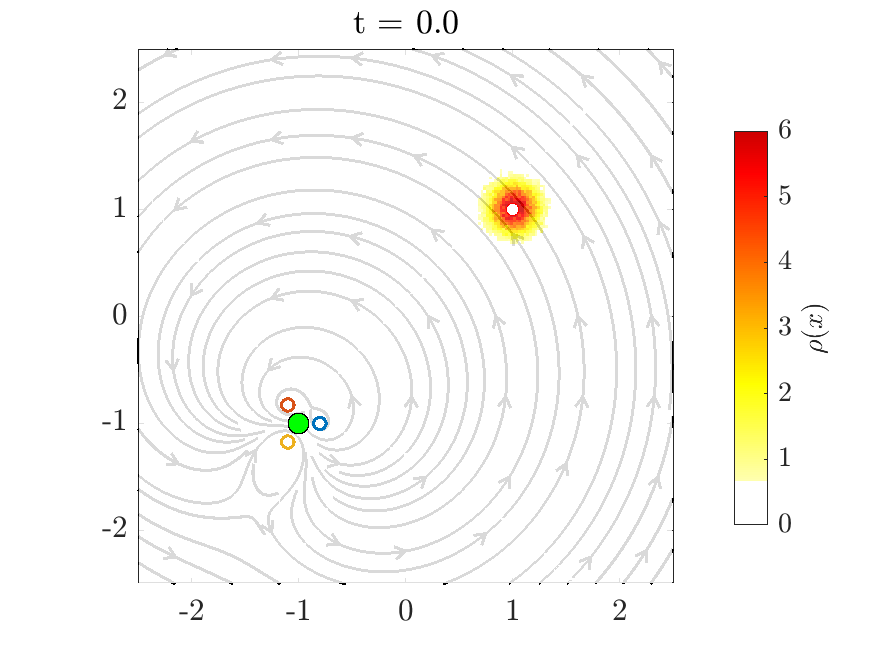}
    \includegraphics[height = 3.0cm,trim={65 30 90 20},clip]{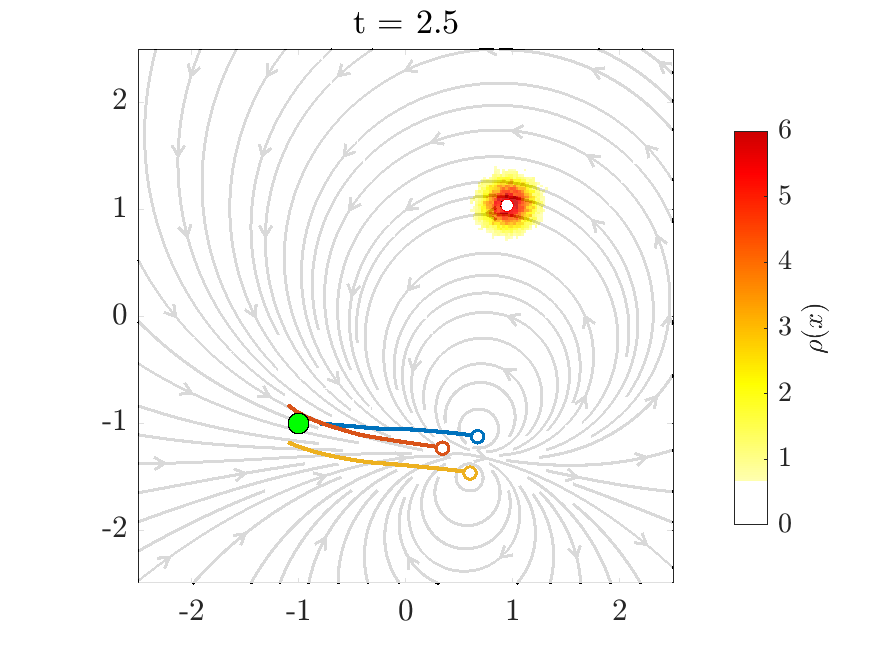}
    \includegraphics[height = 3.0cm,trim={65 30 90 20},clip]{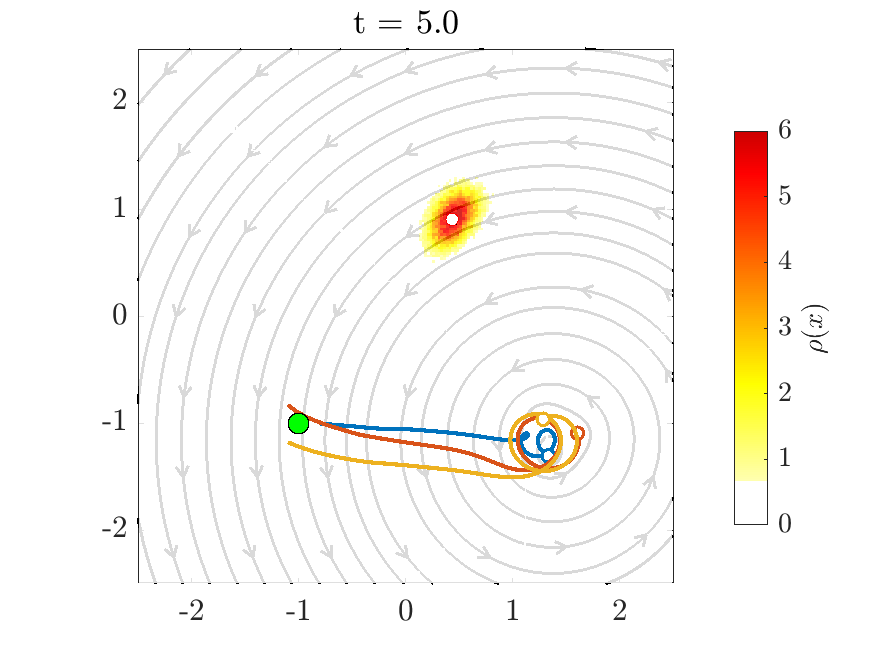}
    \includegraphics[height = 3.0cm,trim={65 30 90 20},clip]{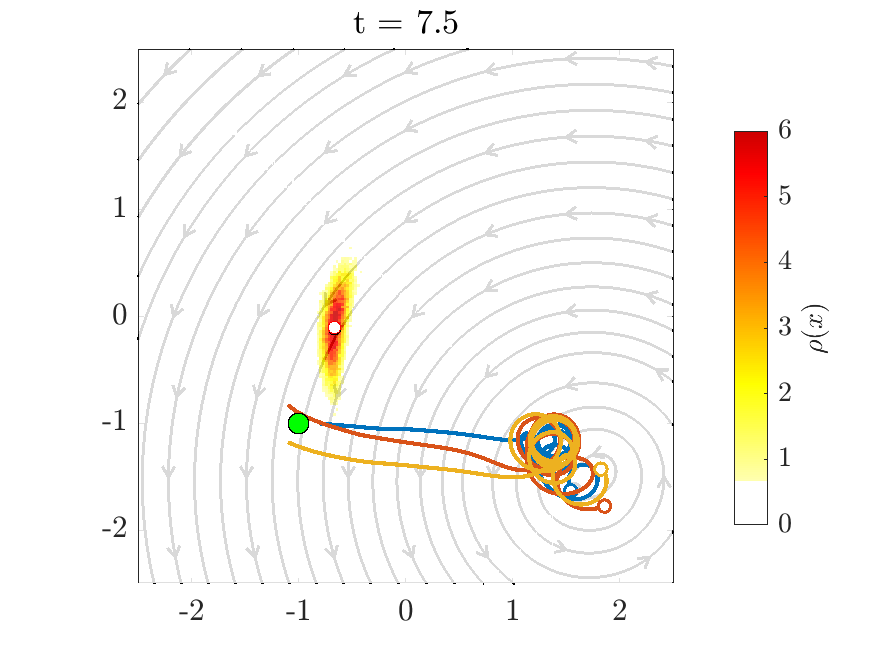}
    \includegraphics[height = 3.0cm,trim={65 30 20 20},clip]{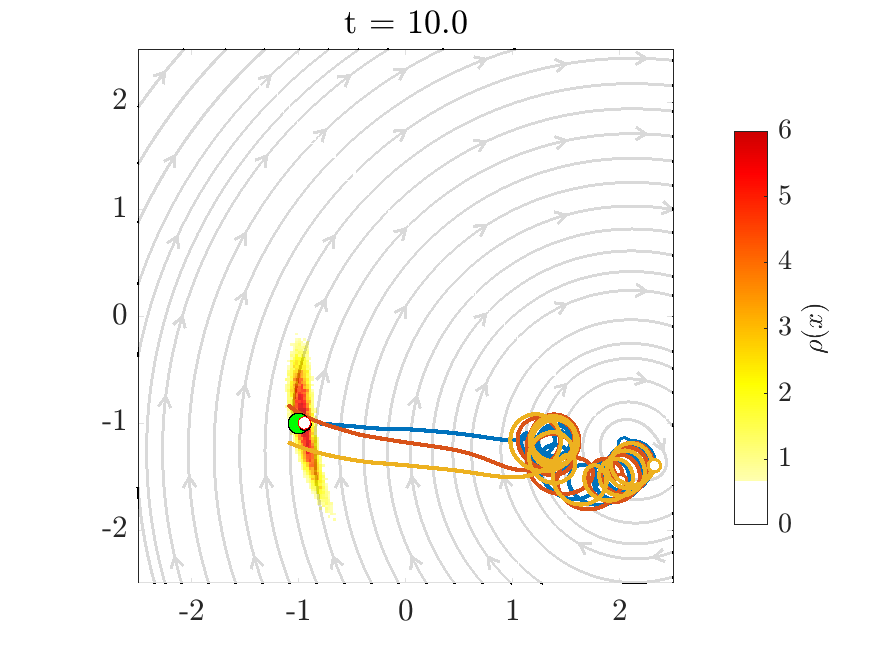}
    \\[3pt]
    \includegraphics[height = 3.0cm,trim={45 30 90 20},clip]{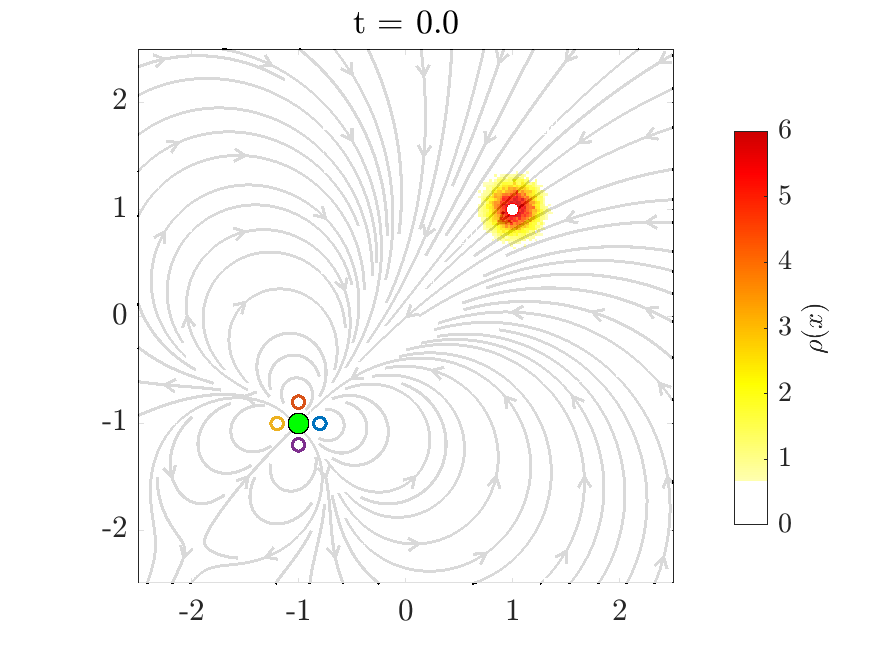}
    \includegraphics[height = 3.0cm,trim={65 30 90 20},clip]{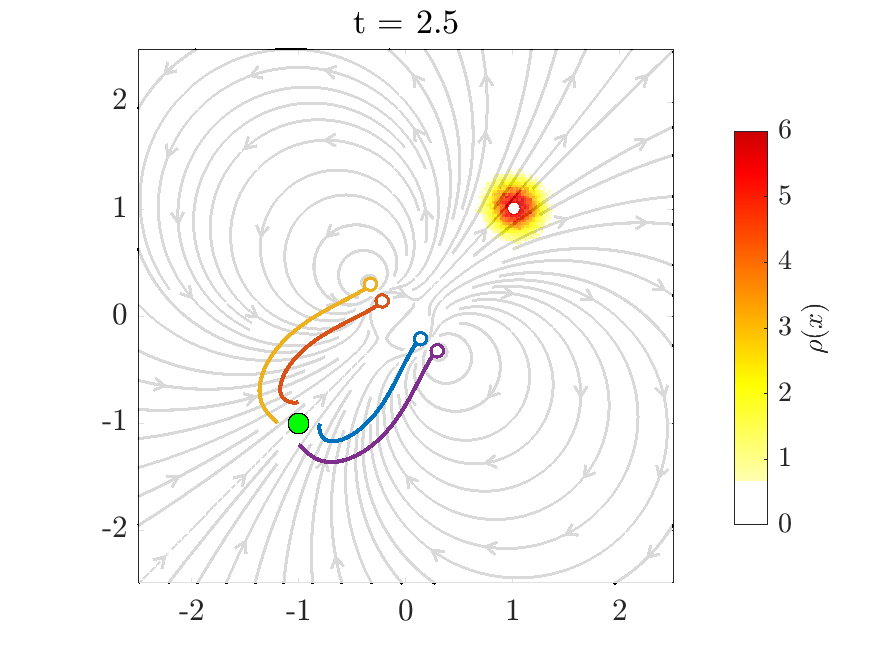}
    \includegraphics[height = 3.0cm,trim={65 30 90 20},clip]{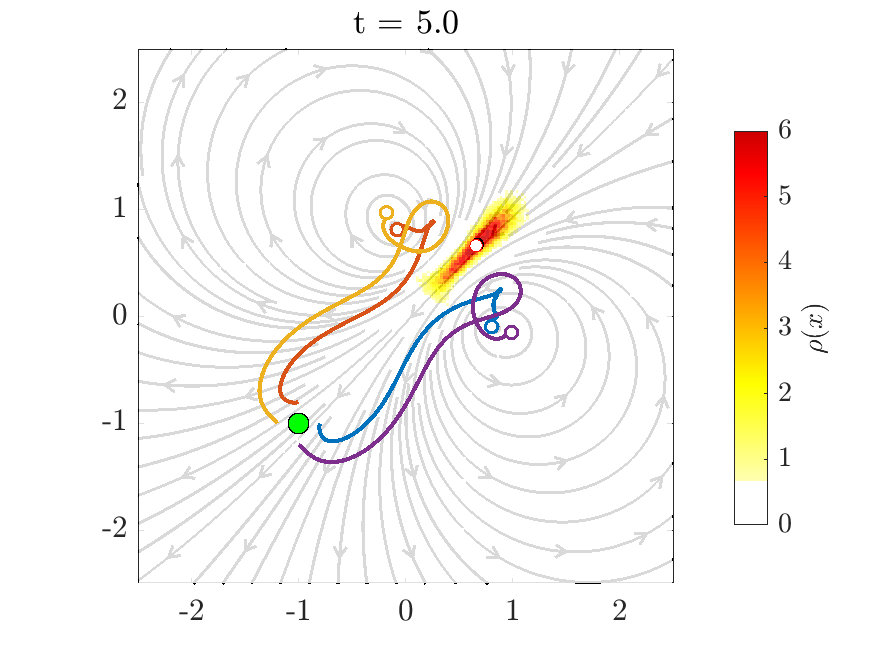}
    \includegraphics[height = 3.0cm,trim={65 30 90 20},clip]{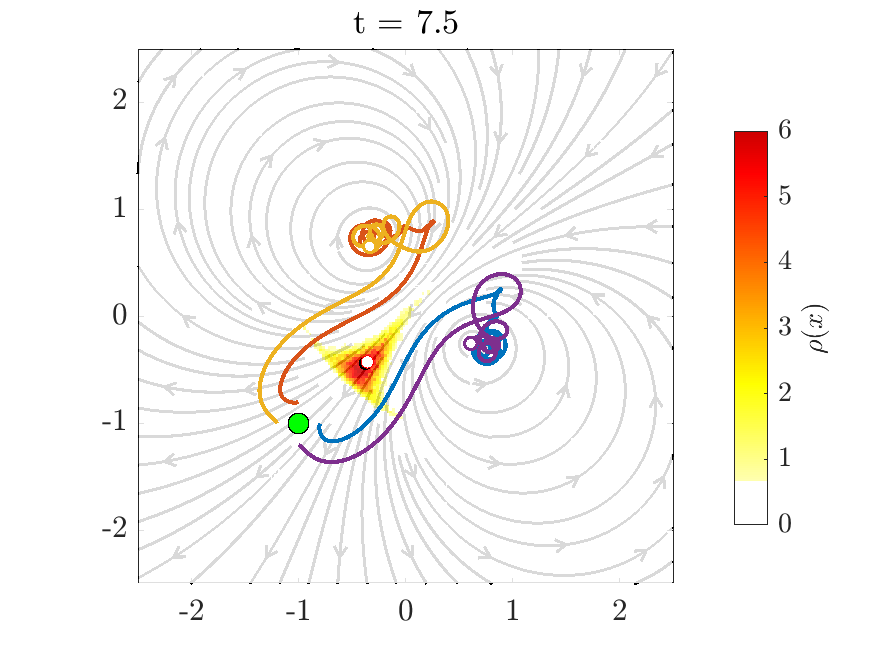}
    \includegraphics[height = 3.0cm,trim={65 30 20 20},clip]{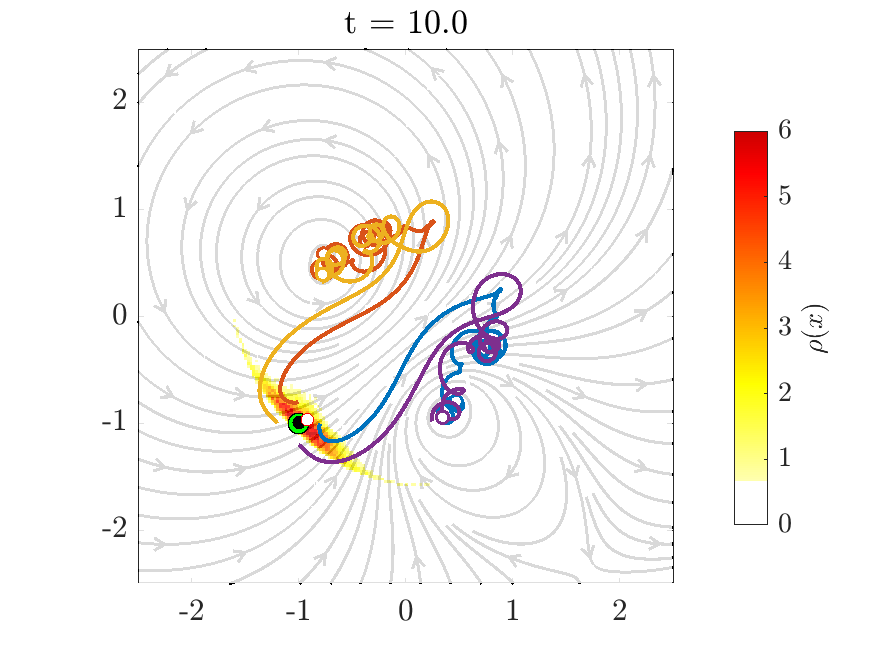}
    \\[3pt]
    \includegraphics[height = 3.34cm,trim={45 0 90 20},clip]{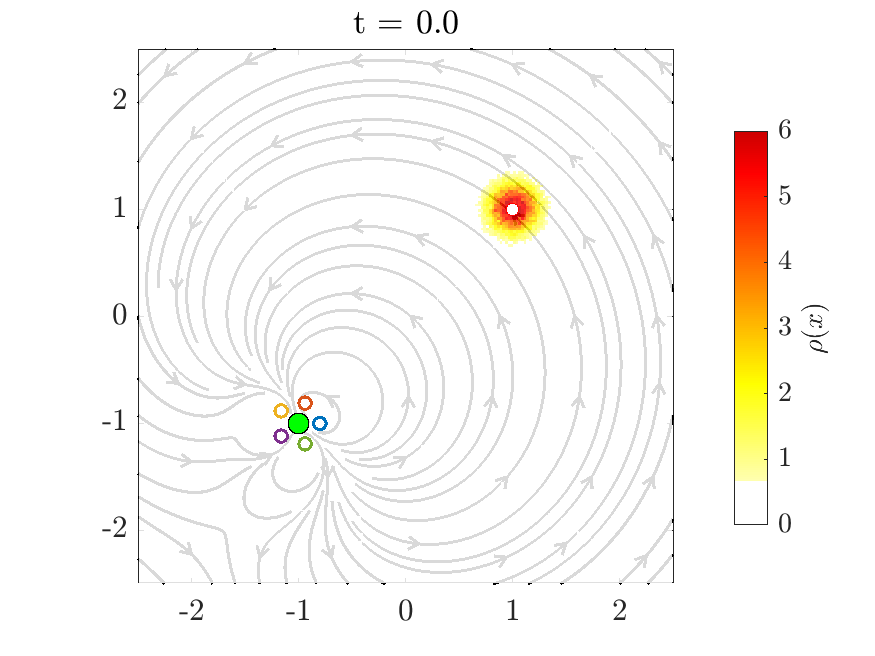}
    \includegraphics[height = 3.34cm,trim={65 0 90 20},clip]{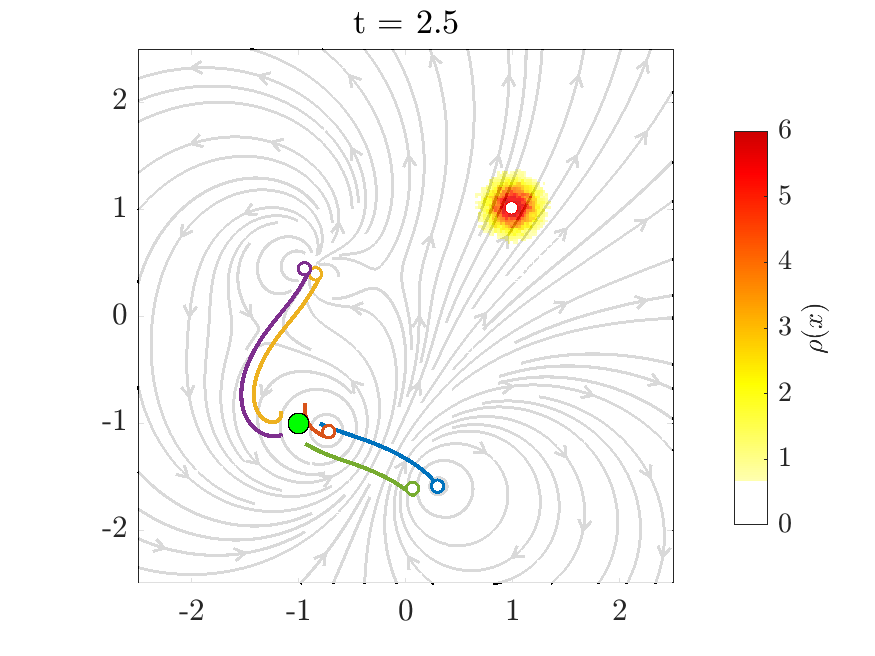}
    \includegraphics[height = 3.34cm,trim={65 0 90 20},clip]{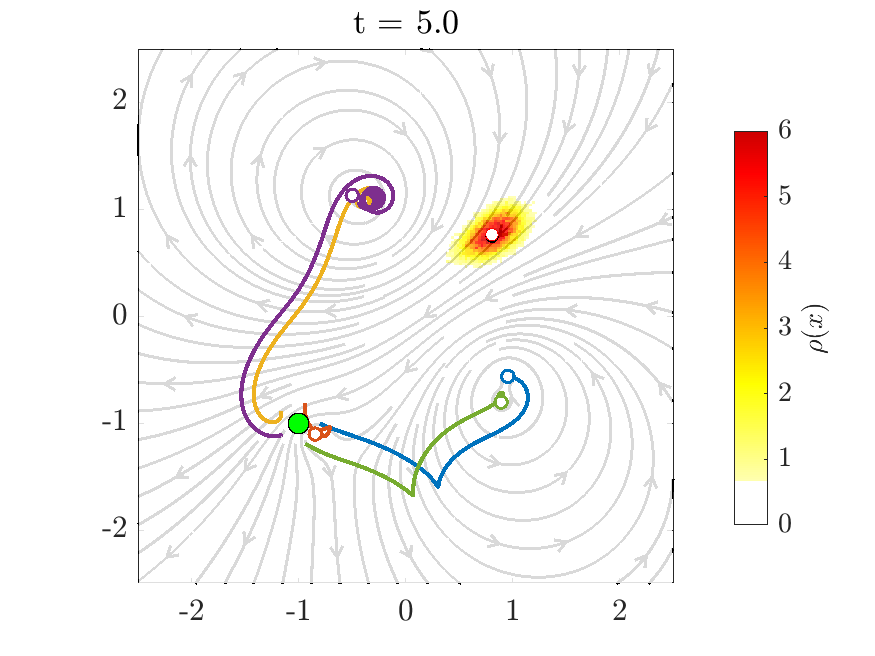}
    \includegraphics[height = 3.34cm,trim={65 0 90 20},clip]{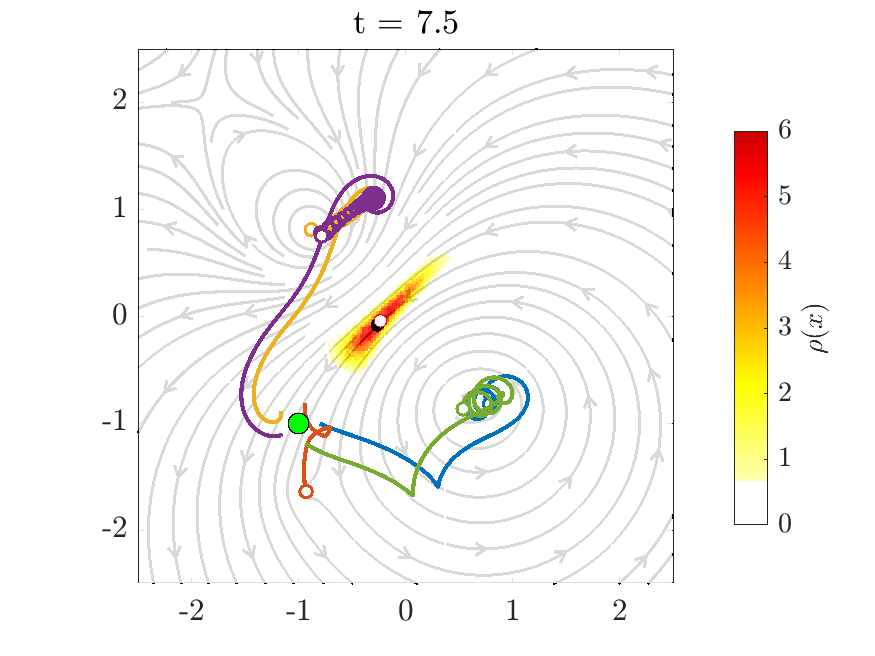}
    \includegraphics[height = 3.34cm,trim={65 0 20 20},clip]{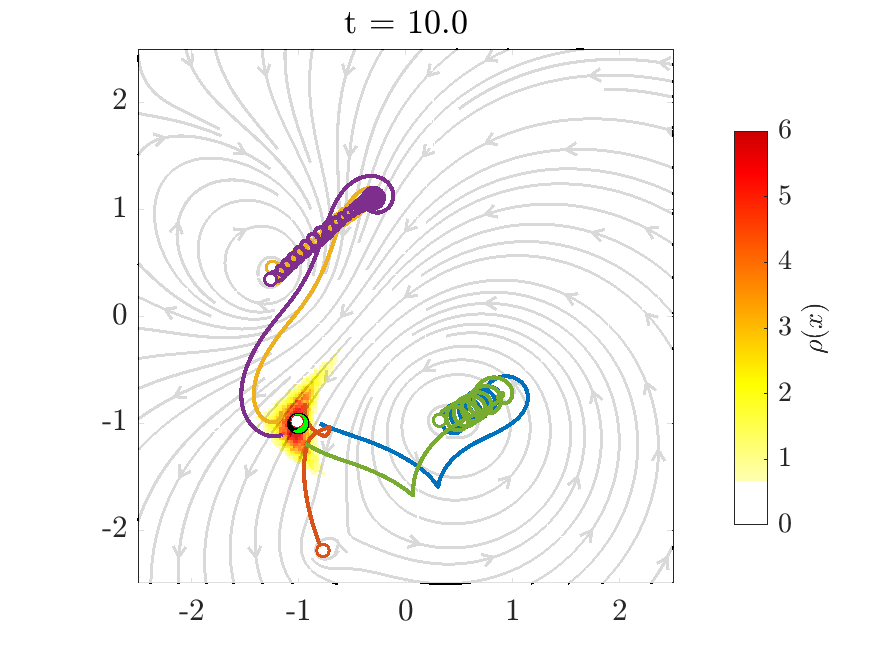}
    \end{center}
    \caption[Optimal solution for torque-only control with varying number of rotors and final time $t_f = 10$ shown as snapshots from the time sequence.]{ Optimal solution for torque-only control with, $n_r=2$ (first row), $n_r = 3$ (second row), $n_r = 4$ (third row), and $n_r = 5$ (fourth row), and final time $t_f = 10$ shown as snapshots from the time sequence.  Streamlines show the direction of the instantaneous velocity fields at each instant shown.  Color shows the fluid particle density at each time instant.  Colored lines show the rotor paths, with circles indicating the rotor position at each instant. 
    The black circle shows the sample mean and the white circle shows the gPC-predicted mean. 
    Each row corresponds to a different number of rotors and each column corresponds to a different time instant, as displayed at the top. }
    \label{fig:lonly_traj}
\end{figure*}

Fig. \ref{fig:lonly_traj} shows the trajectories resulting from the optimal control as computed using DDP for each of these cases and Fig. \ref{fig:lonly_control} shows the corresponding control sequences.   
In the case of $n_r=2$, we see that the rotors translate together by applying strengths of opposite sign in the early part of the trajectory.  Once the rotors have moved closer to the particle distribution, they begin to apply strengths of the same sign, which leads the rotors to spiral about one another, producing a velocity field effectively similar to that of a single rotor in the far field.  In this way, the solution is qualitatively similar to the the solution for the $n_r = 1$ case where the rotor velocity was directly controlled.  In the case where $n_r = 3$ here, the rotors follow a similar pattern, first translating together before spiralling about one another to produce a velocity field similar to that of a single rotor.  In the case of $n_r = 4$, we see a different behaviour, in which the group of $4$ rotors splits into two groups of two, with one pair staying on each side of the line between the target point and the initial distribution location.  If each of these rotor pairs is thought of as a single rotor, the solution trajectory is qualitatively similar to the case of $n_r = 2$ where the rotor velocities were controlled directly.  The two pairs of rotors act symmetrically to retrieve the distribution and push it toward the target. 
This pattern continues for the case of $n_r = 5$, where the additional rotor breaks away from the two pairs, quite similarly to the case of $n_r = 3$ rotors from the previous section, where the last rotor stays near to the target point and exhibits a torque to improve the shape of the distribution as it converges to the target.

\clearpage
\section{Conclusion}
In this work, we have demonstrated a method for numerically computing an optimal control to transport a distribution of fluid particles from an initial distribution to a desired target distribution, while minimizing the variance to reduce stretching, folding, and mixing effects.  
Our method borrows ideas from uncertainty quantification, which have recently been used in a stochastic control context, in order to model and control the propagation of the distribution.  
We apply this method to control a distribution of fluid particles using a system of microrotors, modelled as rotlet singularities, which are allowed to move.  We consider two models of the rotor motion, first considering a case where the rotor translational velocities can be prescribed directly, and second, a case where the rotors are constrained to move according to the flow field produced by the other rotors.  The latter case is significantly more challenging, as the rotor motion is given by a highly coupled, nonlinear system with less actuation than in the former case.  We show that the polynomial chaos with  differential dynamic programming approach allows us to solve these problems numerically to determine a control to effectively transport the distribution.  We analyze the solutions in terms of the optimal cost as a function of the number of rotors and the allotted time, showing that beyond $n_r=4$ rotors and beyond a final time $t_f = 8$, there is negligible improvement to the cost with increasing time or number of rotors.  Further, we analyze the structure of the time-varying flow produced by the optimal control by computing Lagrangian coherent structures of this flow.  Through this analysis, we showed that the optimal control generates a transport barrier which entrains a region of the phase space enclosing the particle distribution which is pulled inward toward the target.  Finally, in solving the case where only the rotor torques are controlled, we show that the optimal control solution tends to exploit pairwise interactions between rotors in which pairs or small groups of rotors interact hydrodynamically to translate together to retrieve the distribution and then spiral about one another to create a rotary flow field which pushes the distribution to the target.  Effectively, these pairs of rotors apply their torques to achieve resultant translation and strength patterns which are qualitatively similar to those of a single rotor in which the velocities can be directly prescribed. 

Exciting future directions for this work include extending these results to the realistic swimmer models as considered in Refs. 
\cite{Buzhardt2019, Buzhardt2020, Buzhardt2021}, where the rotor motion and interactions involve rigid body dynamics and fluid interactions in three dimensions or interactions with boundaries, as considered in \cite{Buzhardt2024}, or where the torques on all rotors are determined by a single, global control input such as a uniform magnetic field.  In these scenarios, the methods developed here could also be used to develop control strategies which explicitly account for other known parametric uncertainties, such as those associated with a swimmer's magnetic moment as considered in \cite{Buzhardt2019}, or uncertainties arising from variations in swimmer geometry which effect the hydrodynamic mobilities.  
In terms of methods and formulation, further developing the connection between the optimal control solution and the Lagrangian coherent structures associated with the optimal flow field remains and interesting line of research.  Such methods could be quite useful for controller design, as this geometric picture gives useful insights into the underlying structure of the time-varying flow induced by the optimal control. 

\appendix

\section{Differential Dynamic Programming} \label{sec:ddp}
Differential dynamic programming computes a locally optimal control around a nominal trajectory by minimizing a quadratic approximation of the value function along this trajectory, and then iteratively optimizing about the new trajectories obtained by applying the locally optimal control.  Here we briefly review the procedure.  First define the value function $V(X_t,t)$ at time $t$ as,
\begin{equation}
V(X_t,t) = \min_{u_t} [ l(X_t,u_t) + V(X_{t+1},t+1)]   
\end{equation}
which expresses the optimal cost-to-go from $X_t$, where the value at the terminal state is defined as $V(X_H,H) = l_H(X_H)$.
Denote by $Q(\delta X, \delta u)$ the change in the value function due to applying change in control input $\delta u$ about the nominal trajectory and consider its quadratic approximation 
\begin{equation}
    Q(\delta X, \delta u) \approx
    \:Q_{X}\delta X + Q_u\delta u + \delta X^T Q_{X u}\delta u + \frac{1}{2}\delta X^TQ_{XX}\delta X + \frac{1}{2}\delta u^T Q_{uu}\delta u 
\end{equation}
where these derivatives are given by 
\begin{subequations}\label{eq:ddp_qgrads}
\begin{align}
Q_{X} &= l_{X}+ F_{X}^T V_{X}'\\
Q_{u} &= l_{u}+ F_{u}^T V_{X}'\\
Q_{XX} &= l_{XX}+ F_{X}^T V_{X X}'F_{X} + V_{X}'\cdot F_{XX}\\
Q_{uu} &= l_{uu}+ F_{u}^T V_{XX}'F_{u} + V_{X}'\cdot F_{uu}\\
Q_{Xu} &= l_{Xu}+ F_{X}^T V_{X X}'F_{u} + V_{X}'\cdot F_{Xu}
\end{align}
\end{subequations}
where the notation $(\cdot)'$ indicates the next time step.
The algorithm proceeds by computing these derivatives by recursing backward in time along the nominal trajectory from the end of the horizon.
At each iteration, the control policy is improved by optimizing this quadratic expansion with respect to $\delta u$
\begin{equation}
    \delta u^* = \arg \min_{\delta u}Q(\delta X,\delta u) = -Q_{uu}^{-1}\left(Q_u + Q_{u X}\delta X\right)
\end{equation}
This can be seen as providing a descent direction in the space of control policies.  An updated nominal control is then computed by a linesearch over a stepsize parameter $\alpha$ to update the policy, that is 
\[
u_{\text{new}} = u - \alpha Q_{uu}^{-1}Q_u - Q_{uu}^{-1}Q_{u X}\delta X
\] 
and this new control is applied to obtain a new nominal trajectory, and this procedure is iterated until the relative change in cost falls to less than a specified tolerance. Full details of the algorithm can be found in Refs. \cite{Tassa2012,Yakowitz1984} or in Ref. \cite{Boutselis2019} for an implementation in terms of gPC dynamics.

\subsection{Jacobians and Hessians of gPC dynamics}

In Eq. \ref{eq:ddp_qgrads}, the Jacobians and Hessians of the polynomial chaos form of the dynamics in Eq. \ref{eq:gpc_ode_compact} are required.
To derive these, it is helpful to first write the  dynamics of Eq. \ref{eq:gpc_ode_compact} more explicitly as  
\[
\frac{d}{dt}
\begin{bmatrix}
    x_{11}\\ 
    x_{12}\\ 
    \vdots \\ 
    x_{1K}\\ 
    x_{21}\\ 
    \vdots\\ 
    x_{2K}\\ 
    \vdots \\ 
    x_{nK}
\end{bmatrix}
=
\begin{bmatrix}
    \<f_1,\phi_1\>/\<\phi_1^2\> \\
    \<f_1,\phi_2\>/\<\phi_2^2\> \\
    \vdots \\ 
    \<f_1,\phi_K\>/\<\phi_K^2\> \\[0.5ex] 
    \<f_2,\phi_1\>/\<\phi_1^2\>  \\ 
    \vdots\\ 
    \<f_2,\phi_K\>/\<\phi_K^2\>\\ 
    \vdots \\ 
    \<f_n,\phi_K\>/\<\phi_K^2\>
\end{bmatrix}
\]
From this, it can be seen that the needed gradients of $\mathbf{f}$ with respect to the gPC state can be expressed in terms of the Jacobians and Hessians of the original system, $f$, as
\begin{align}
\frac{\p \mathbf{f}_l}{\p x_{ab}} 
&= \frac{\<\frac{\p f_i}{\p x_{ab}}\phi_j\>}{\<\phi_j^2\>}
= \frac{\<\frac{\p f_i}{\p x_{a}}\phi_b\phi_j\>}{\<\phi_j^2\>} \\[1ex]
\frac{\p^2 \mathbf{f}_l}{\p x_{ab}\p{x_{cd}}} &= \frac{\<\frac{\p^2 f_i}{\p x_{ab}\p x_{cd}}\phi_j\>}{\<\phi_j^2\>}
= \frac{\<\frac{\p^2 f_i}{\p x_{a}\p x_{c}}\phi_b\phi_d\phi_j\>}{\<\phi_j^2\>} 
\end{align}
where $i,a,c=1,\dots,n$, and $j,b,d=1,\dots,K$, and $l=K(i-1)+j$ is an index ranging from $1,\dots,nK$. 
Similarly, the control gradients are 
\begin{subequations}
\begin{align}
    \frac{\p \mathbf{f}_l}{\p u_g} &=  
  \frac{\<\frac{\p f_i}{\p u_g}\phi_j\>}{\<\phi_j^2\>} 
  \\
  \frac{\p^2 \mathbf{f}_l}{\p u_g \p u_h} &=  
  \frac{\<\frac{\p^2 f_i}{\p u_g \p u_h}\phi_j\>}{\<\phi_j^2\>} 
  \\
  \frac{\p^2\mathbf{f}_l}{\p x_{ab} \p u_g} &= \frac{\<\frac{\p^2 f_i}{\p x_{a}\p u_g}\phi_b\phi_j\>}{\<\phi_j^2\>}
  \end{align}
\end{subequations}
where $g,h = 1, \dots, n_c$. 

In Eq. \ref{eq:ddp_qgrads} the gradients of the loss function with respect to the state and control are also needed.  
For the quadratic output tracking loss function in Eq. \ref{eq:oc_loss}
\begin{equation}
    l(X_t,u_t) = (M(X_t)-\yr_t)^TS(M(X_t)-\yr_t) + u_t^TRu_t
\end{equation}
where $\yr$ is a vector of the relevant target moments to be tracked. 
The gradients of this needed for DDP are as follows (assuming that the weighting matrix $S$ is diagonal).
\begin{subequations}
\begin{align}
    \frac{\partial l}{\p X_i} &= 2S_{jj}\frac{dM_j}{dX_i} - 2S_{jj}\yr_j\frac{dM_j}{dX_i} \label{eq:lX}\\
    \frac{\p^2 l}{\p X_i \p X_j} &= 2S_{kk}\left(\frac{\p^2M_k}{\p X_i \p X_j} M_k + \frac{\p M_k}{\p X_i}\frac{\p M_k}{\p X_j}  \right)  \label{eq:lXX}\\
    \frac{\p^2 l}{\p X \p u} &= 0 \\
    \frac{\p l}{\p u} &= 2Ru\\
    \frac{\p^2 l}{\p^2 u} &= 2R
\end{align}
\end{subequations}
where summation over repeated indices is implied in Eq. \ref{eq:lX} and \ref{eq:lXX}.


 \bibliographystyle{elsarticle-num} 
 \bibliography{gpc_rotlets}


\end{document}